\documentclass[bezier,a4paper,leqno,draft]{amsart}
\usepackage[latin1]{inputenc}
\usepackage{amsmath}
\usepackage{latexsym}
\usepackage{amsthm}
\input{psfig}
\newcommand{\bmath}[1]{\mbox{\mathversion{bold}$#1$}}

\newcommand{\bfC}{\bmath{C}}

\newcommand{\R}{\bmath{R}}
\newcommand{\bfR}{\bmath{R}}
\newcommand{\bfH}{\bmath{H}}

   \newtheorem{theorem}{Theorem}
   \newtheorem{lemma}[theorem]{Lemma}
   \newtheorem{definition}[theorem]{Definition}
\title[Minimal surfaces in $\bfR^3$ vs. CMC $1$ surfaces in 
$\bfH^3$]{Mean curvature 
one surfaces in hyperbolic space, and their relationship to minimal 
surfaces in Euclidean space
}
\date{\today}
\author{Wayne Rossman}
\address{%
   Department of Mathematics, Faculty of Science,
   Kobe University,
   Rokko, Kobe 657-8501, Japan 
}
\begin{document}
\maketitle

\begin{abstract}
We describe local similarities and global 
differences between minimal surfaces in Euclidean $3$-space and 
constant mean curvature $1$ surfaces in hyperbolic $3$-space.  
We also describe how to solve global period problems for constant mean 
curvature $1$ surfaces in hyperbolic $3$-space, and we give an overview 
of recent results on these surfaces.  We include computer graphics of 
a number of examples.  
\end{abstract}

\section*{Introduction}

In this article we explain some aspects of research on constant mean 
curvature $1$ surfaces in hyperbolic $3$-space.  The first half of the 
article emphasizes the most elementary elements of the theory, and the last 
half attempts to explain less elementary elements as 
elementarily as possible.  Comparison with minimal surfaces in Euclidean 
$3$-space is used, and such a comparison is 
natural (in fact, almost unavoidable), as ideas in the theory have 
evolved directly from ideas in minimal 
surface theory.  We begin with this comparison in this introduction.  

Minimal surfaces in Euclidean $3$-space $\bfR^3$ have been a subject 
of investigation since the 1800's, as evidenced by the names of nineteenth 
century mathematicians associated with some of the most famous 
complete (i.e. without boundary) examples and with
the subject's primary theorem: Enneper's surface, 
Riemann's singly-periodic staircase surface, Scherk's 
periodic surfaces, and the Weierstrass representation theorem.  
The first half of the twentieth century, however, saw the focus turn to 
compact minimal surfaces with boundary, in particular, to the Plateau 
problem of determining the compact minimal surfaces with 
given boundary in $\bfR^3$, as seen in the works of Douglas, Rado, and 
Courant.  (Although, in the early 1900's, Bernstein did prove 
that a complete minimal graph must be a plane.)  

But recently there has been much productive 
research on complete minimal surfaces, with the work of Meeks, 
Hoffman, Karcher, Osserman, Ros, Rosenberg, and others.  And the research 
has been taking advantage of a modern-day tool -- computer graphics.  
A striking example of how computer 
graphics has helped in the study of minimal surfaces is 
the Costa surface, whose proof of embeddedness \cite{HoMe} came as a direct 
result of computer graphics experiments (providing a counterexample to 
a longstanding conjecture that the plane, 
catenoid, and helicoid are the only complete embedded 
minimal surfaces with finite topology).  And the recent proofs of 
existence of numerous examples of minimal surfaces is due in part 
to the ability to first check on a computer if existence is plausible.  
(The beauty of the graphics, as well, cannot be ruled out as a 
motivating factor.)  Then, the many examples have given us a better 
understanding of the nature of minimal surfaces.  

Minimal surfaces, and constant mean curvature (CMC) surfaces as well, 
are natural objects of study, as they are critical for area with respect to 
volume preserving variations that fix their boundaries.  And minimal surfaces 
are especially important in the class of CMC surfaces, 
as they are area critical for all variations, not just volume 
preserving ones.  The particular interest in minimal surfaces 
comes also from the elegance of the Weierstrass representation, which 
describes minimal surfaces in terms of a pair of holomorphic functions.  
Actually, nonminimal CMC surfaces can be described with holomorphic 
functions as well \cite{DPW}, but the descriptions are more complicated, so 
minimal surfaces are unique amongst all CMC surfaces in that they have a 
"simple" holomorphic description.  

Recently, the field has also grown in breadth, by influencing related areas.  
As one among many 
examples, it has influenced the study of CMC surfaces in other 
$3$ dimensional space forms (such as the $3$ dimensional hyperbolic and 
spherical spaces $\bfH^3$ and $S^3$).  Regarding this, we 
should first remark on an elementary correspondence between 
CMC surfaces in the different $3$ dimensional space forms, described by 
Lawson \cite{L} and often attributed to him.  
It gives, in particular, a local correspondence between minimal 
surfaces in $\bfR^3$ and CMC $1$ surfaces in $\bfH^3$.  So it is not 
surprizing that there is a kind of "Weierstrass" representation in terms of 
a pair of holomorphic functions for CMC $1$ surfaces in $\bfH^3$ as well, 
which we call the Bryant representation here, after its discoverer \cite{B}.  
Furthermore, 
among all CMC surfaces in $\bfH^3$, CMC $1$ surfaces are unique in having 
a relatively simple description in terms of a pair of holomorphic functions, 
as we expected, since minimal surfaces in $\bfR^3$ are 
similarly endowed.  (However, since CMC $1$ surfaces in $\bfH^3$ are not 
minimal, they are area critical only for volume preserving 
variations -- not all of the variational 
properties of minimal surfaces in $\bfR^3$ are preserved by the Lawson 
correspondence).  

\begin{figure}
\vspace{0.5in}
\centerline{
        \hbox{
		\psfig{figure=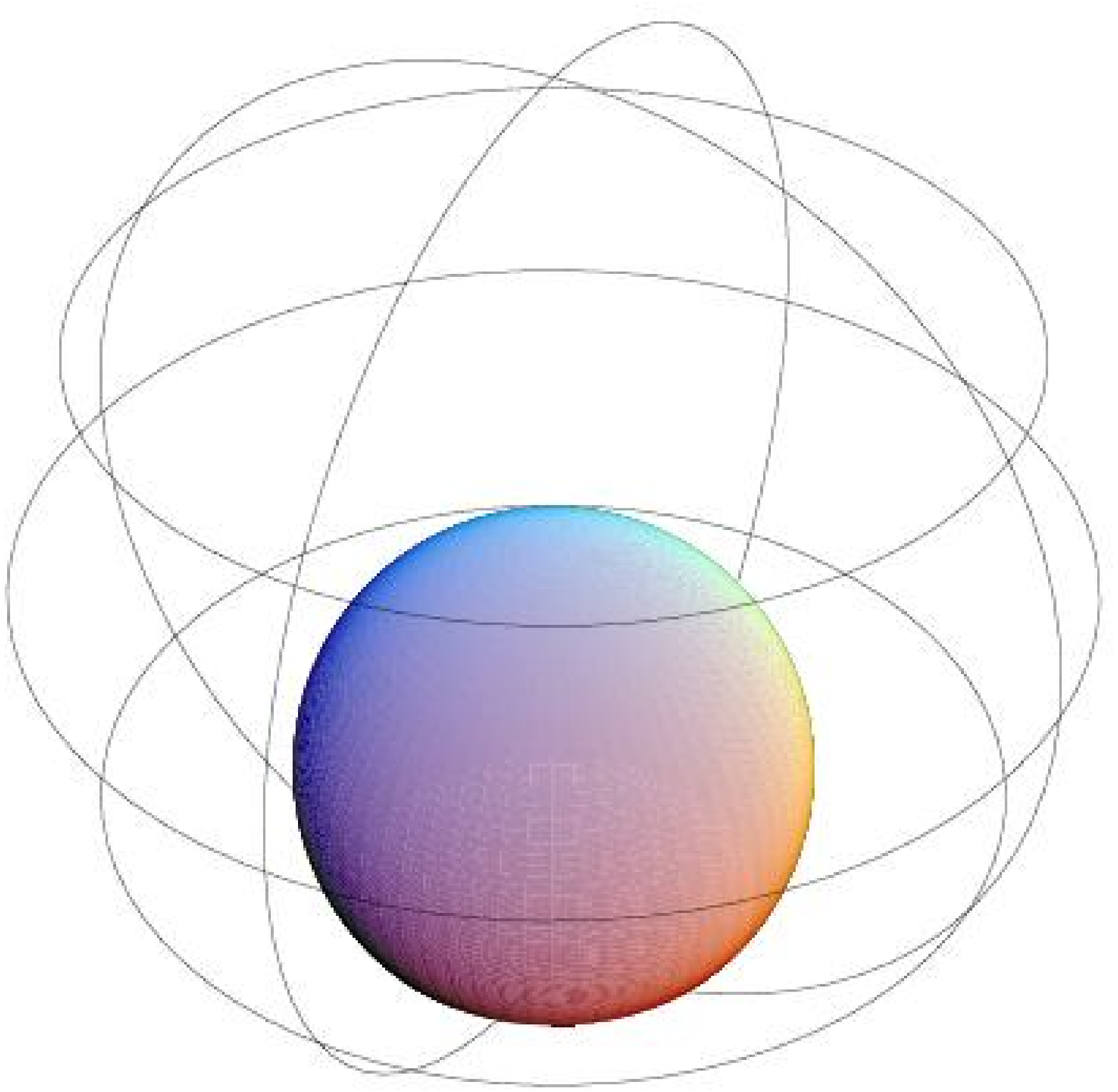,width=2in}
		\hspace{0.5in}
		\psfig{figure=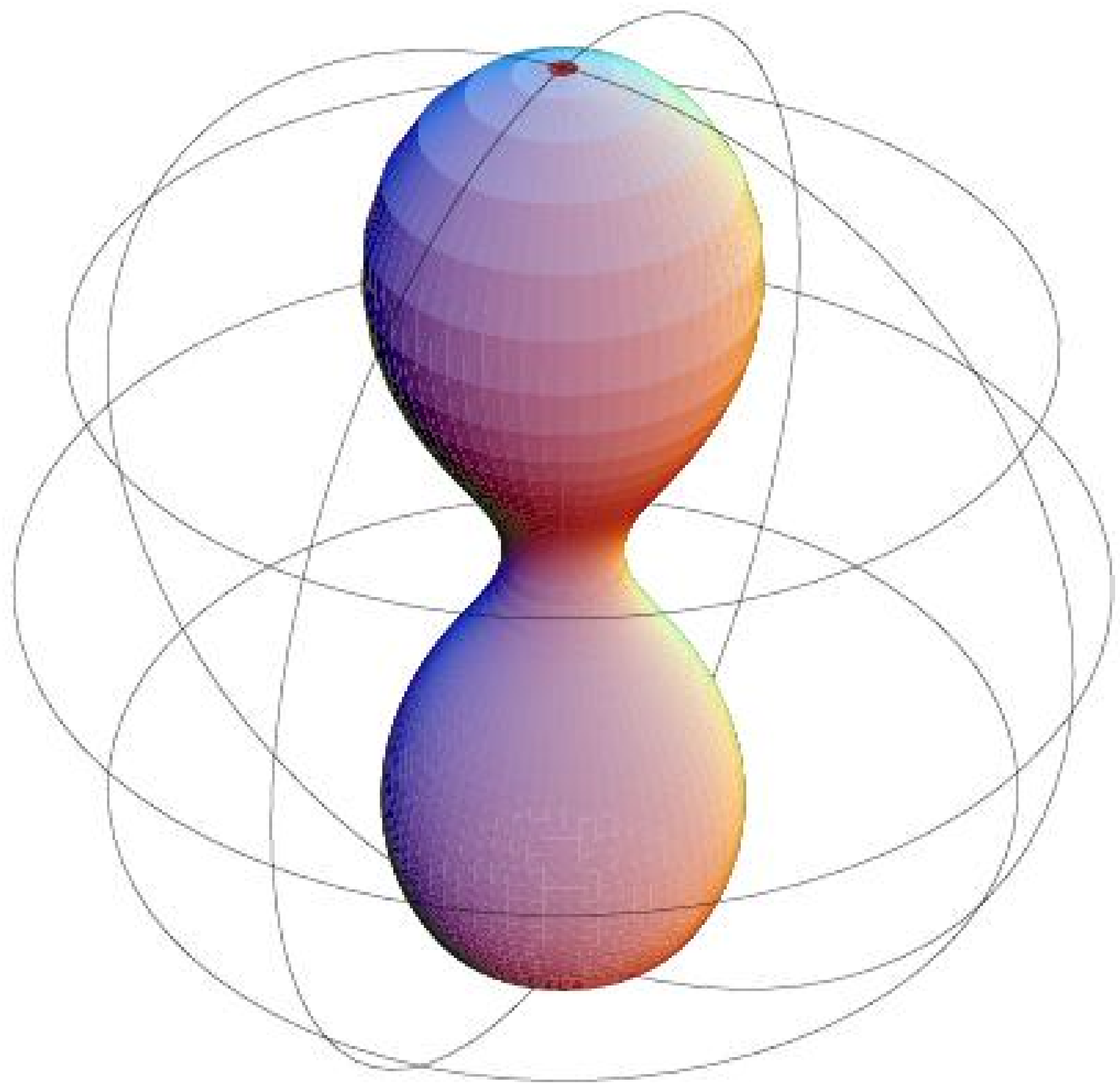,width=2in}
	}
}
\caption{CMC $1$ horosphere and catenoid cousin in $\bfH^3$.  The catenoid 
cousin was first described in \cite{B}.}
\label{figure8}
\end{figure}

\begin{figure}
\vspace{0.5in}
\centerline{
        \hbox{
		\psfig{figure=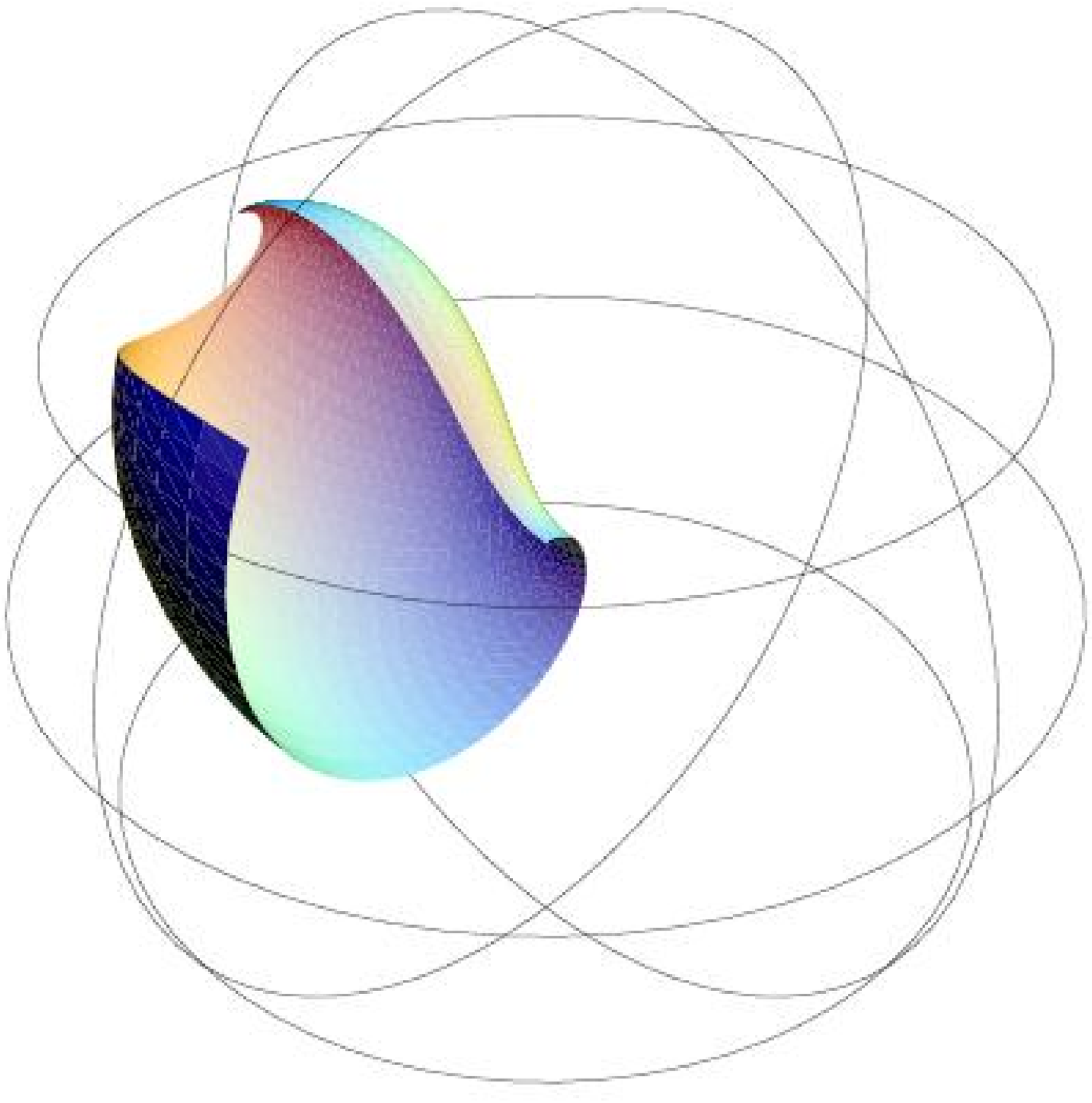,width=2in}
		\hspace{0.5in}
		\psfig{figure=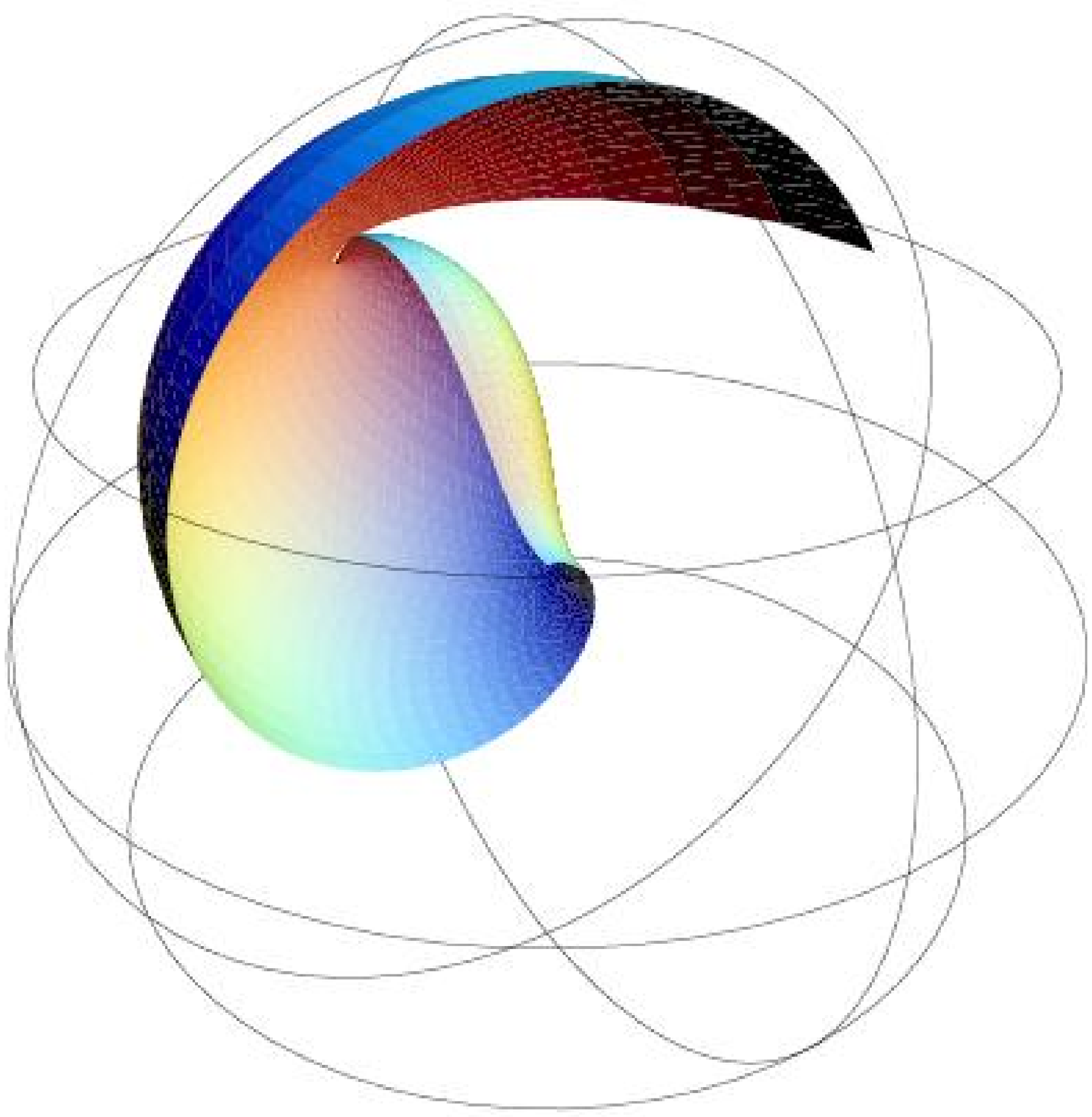,width=2in}}}
		\vspace{0.5in}
		\centerline{
		\hbox{  \psfig{figure=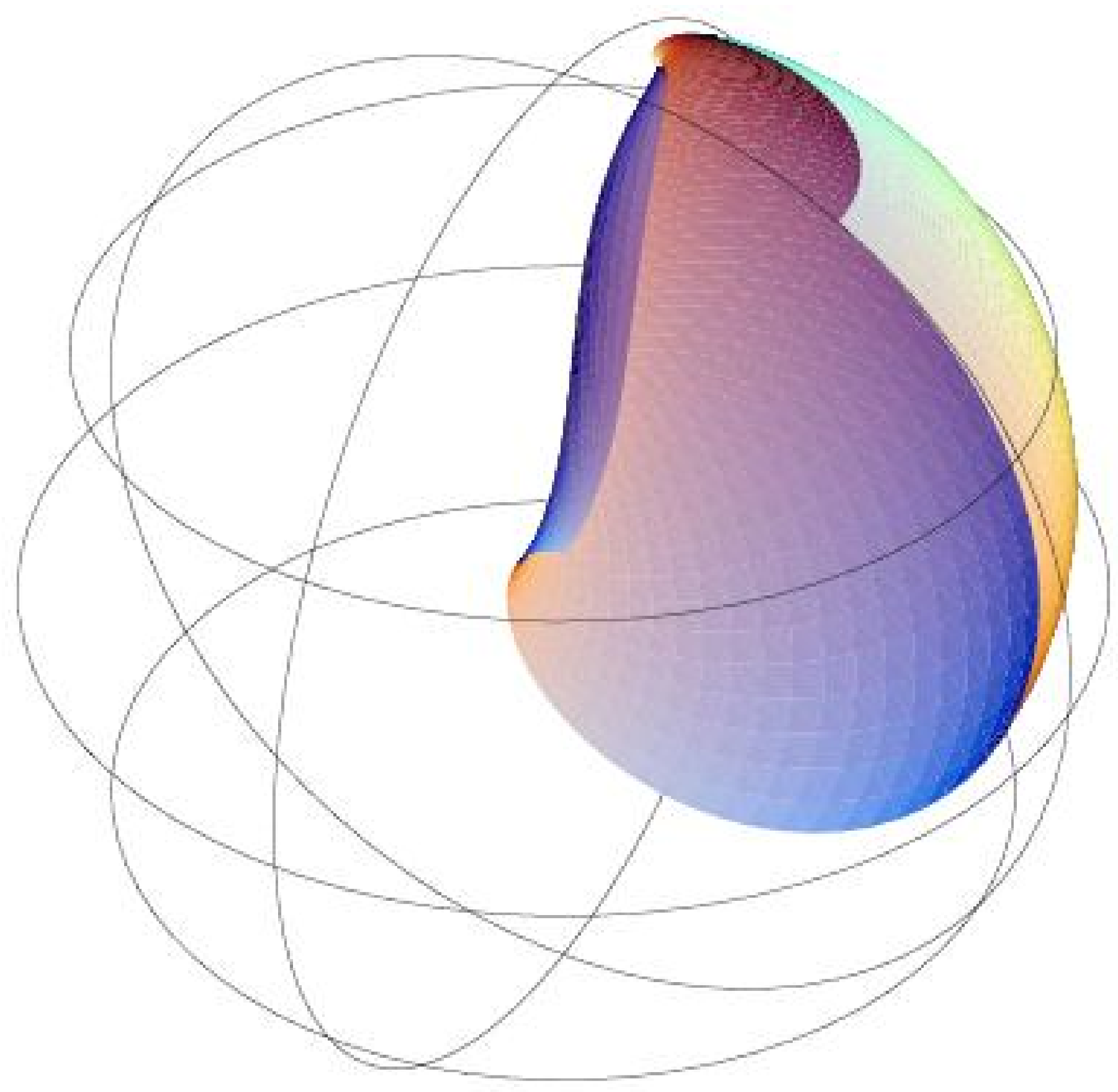,width=2in}
		\hspace{0.5in}
		\psfig{figure=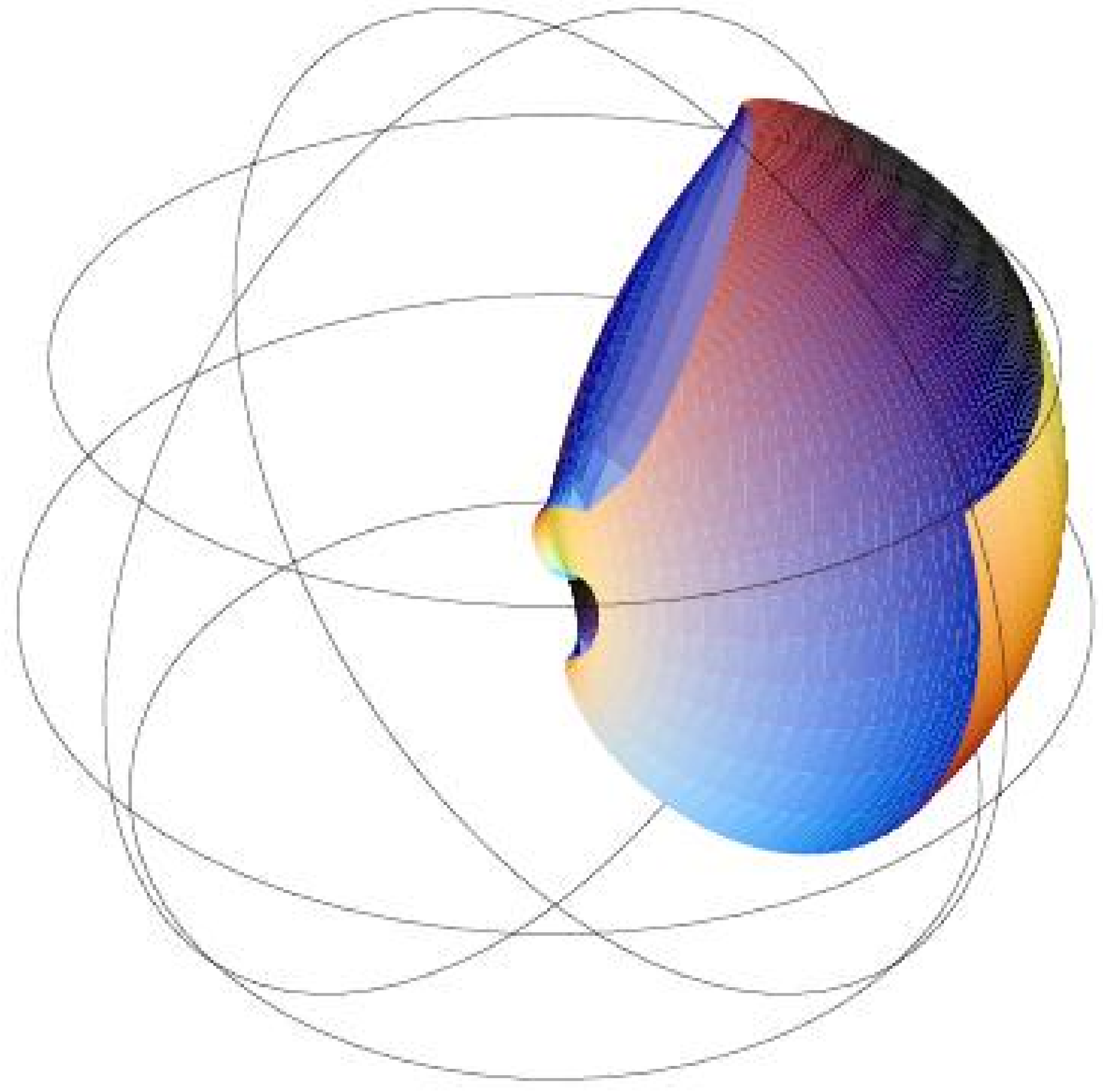,width=2in}
	}
}
\caption{CMC $1$ Enneper cousin \cite{B} (here we show two pictures, the 
second including more of the end) and Enneper cousin dual \cite{ruy1} and 
genus $1$ Enneper cousin dual in $\bfH^3$.  Recently, it has been 
established \cite{ruy3} that the genus $1$ minimal Enneper surface is 
symmetric and nondegenerate, so Theorem 5 implies existence of the third 
surface pictured here.  (Note that only one of four congruent pieces of 
each of these surfaces is shown here, and that the full surfaces can 
be made by reflecting these pieces across their boundary geodesics.)}
\label{figure12}
\end{figure}

\begin{figure}
\vspace{0.5in}
\centerline{
        \hbox{
		\psfig{figure=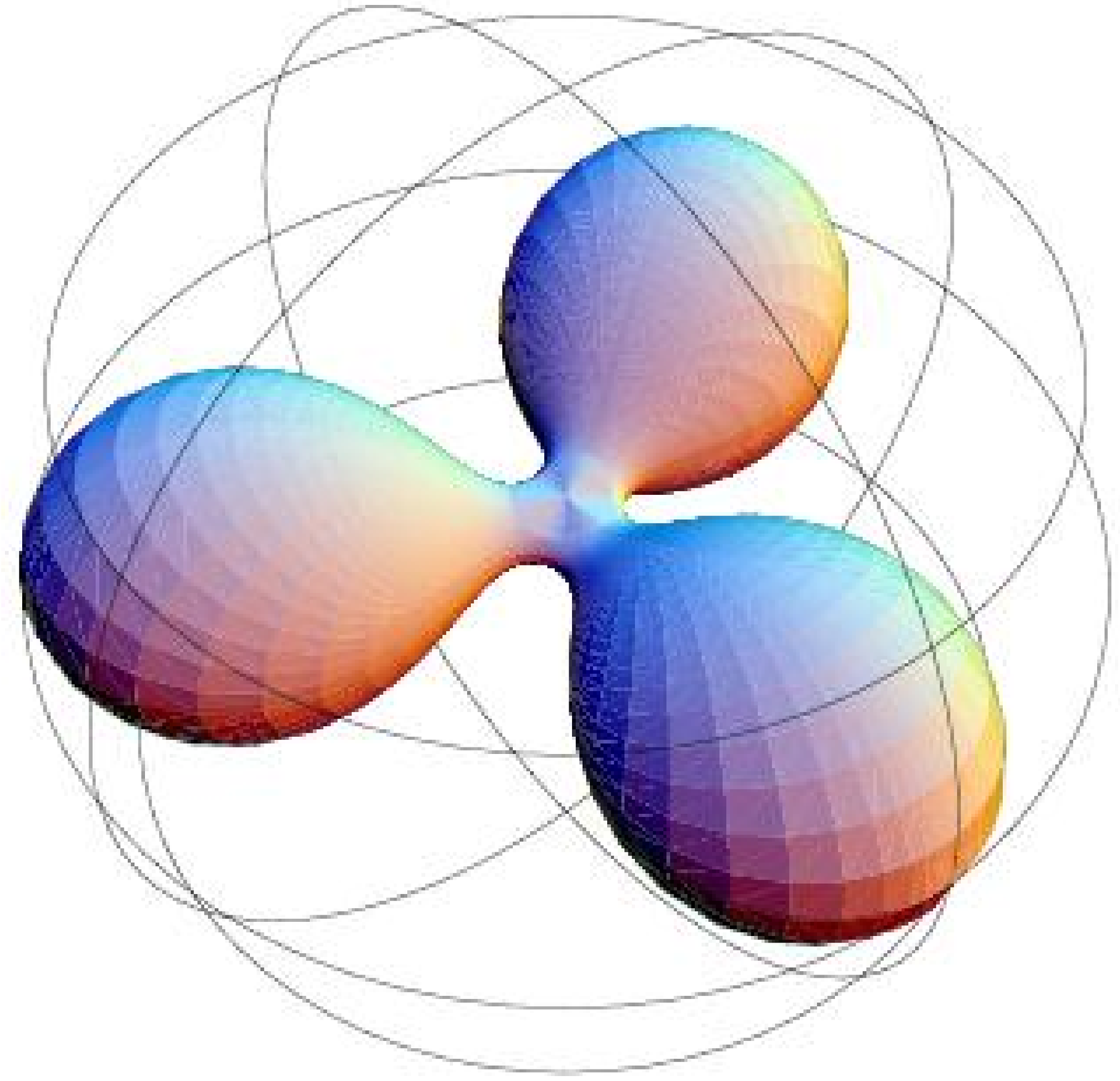,width=2in}
		\hspace{0.5in}
		\psfig{figure=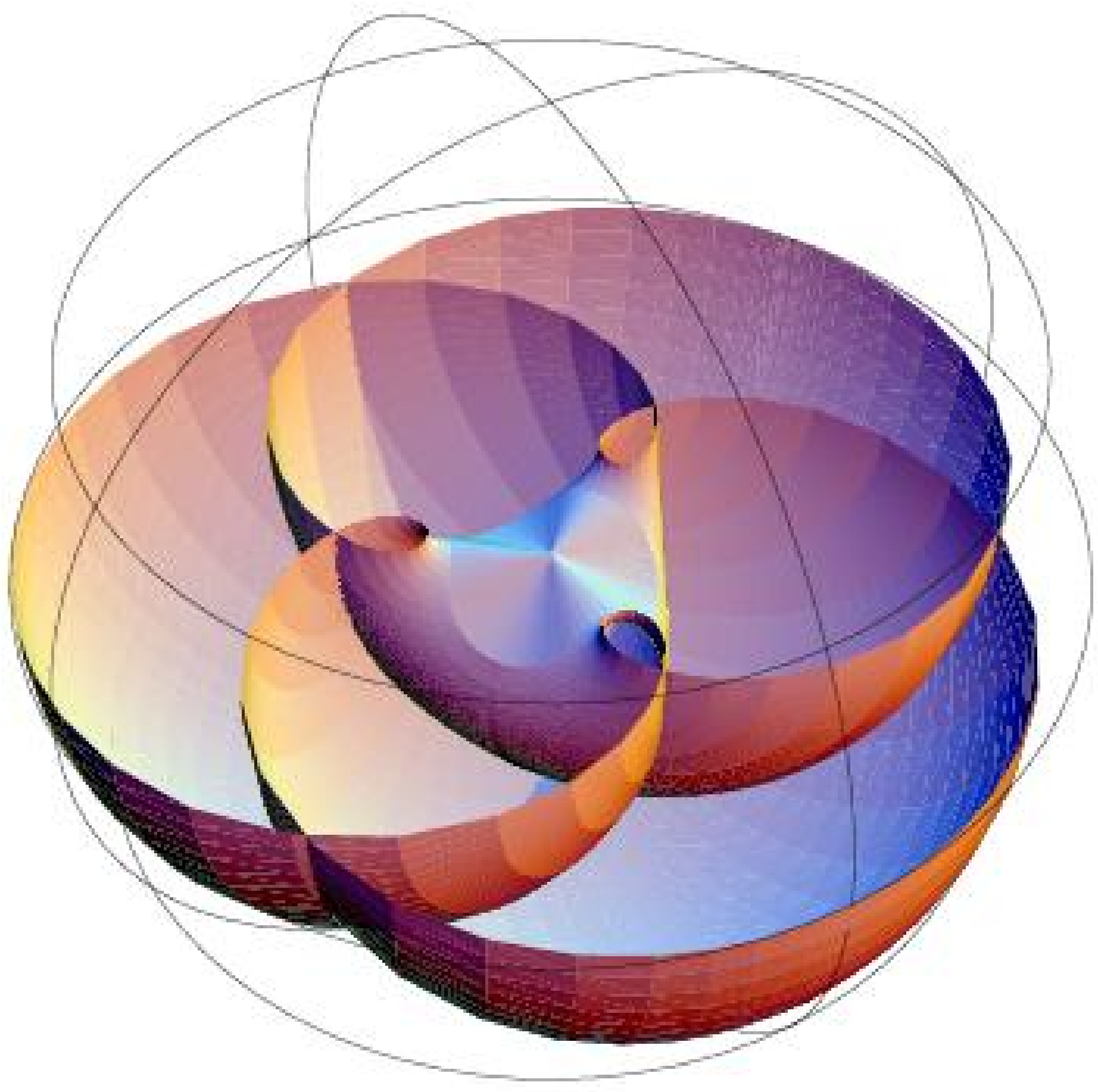,width=2in}}}
		\centerline{\hbox{
		\psfig{figure=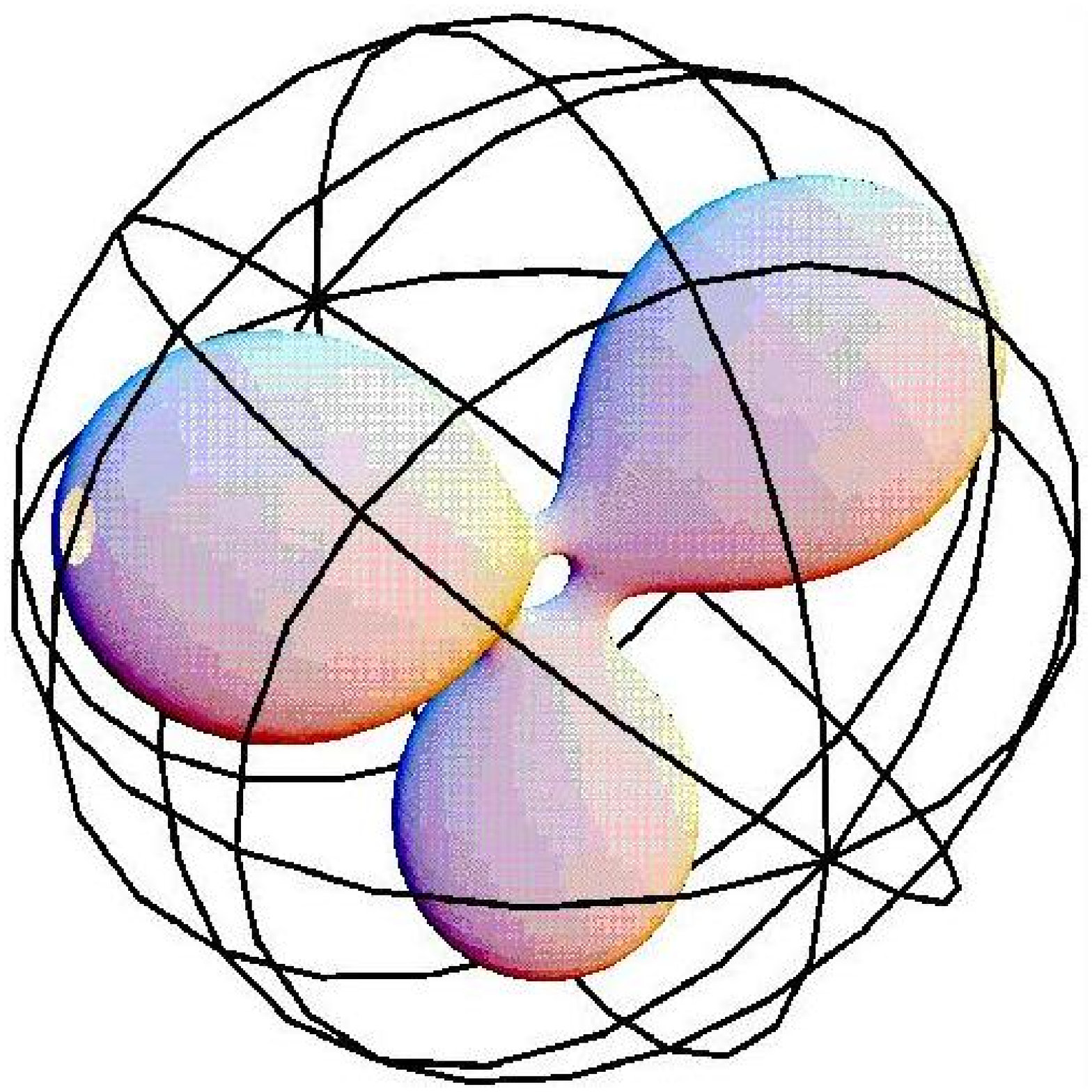,width=2in}
	}
}
\caption{Two different CMC $1$ trinoid cousin duals (proven to exist in 
\cite{uy3}) and a 
CMC $1$ genus $1$ trinoid cousin dual in $\bfH^3$ (for which existence 
follows from Theorem 5).  The graphics for the 
genus $1$ surface here were made by Katsunori Sato of Tokyo Institute 
of Technology.  Although these surfaces are proven to exist, and numerical 
experiments show that some of them are embedded (as two of the 
pictures here are), none have yet been proven to be embedded.}
\label{figure14}
\end{figure}

So locally minimal surfaces in $\bfR^3$ are indeed equivalent to 
CMC $1$ surfaces in $\bfH^3$, and share many analogous properties.  But 
interesting and essential differences arise when considering 
the two types of surfaces globally.  Both the Weierstrass 
and Bryant representations require an integration to produce 
the actual immersions into $\bfR^3$ and $\bfH^3$, which locally can 
always be done, but when we attempt to do it globally on a 
non-simply-connected region, there are period problems, i.e. upon 
integration, we might find that the immersion is well defined only on the 
region's universal cover.  Solvability of 
these period problems (adjusting free parameters to get the immersion 
well defined on the region itself) is a global question, and 
the answer can be different in the two cases.  As a general (but inexact) 
rule of thumb, solvability is more likely in the $\bfH^3$ case, leading 
to a wider variety of surfaces in the $\bfH^3$ case.  
For example, a genus $1$ surface with finite total curvature and two 
embedded ends cannot exist as a minimal surface in $\bfR^3$, but it 
does exist as a CMC $1$ surface in $\bfH^3$.  
And a genus $0$ surface with finite total curvature and two embedded ends 
exists as a minimal surface in $\bfR^3$ only if it is a surface of 
revolution, but it may exist as a CMC $1$ surface in $\bfH^3$ without being 
a surface of revolution.  

Another difference between these two types of surfaces is that there 
are no relations regarding embeddedness of the surfaces in the two different 
cases.  For example, the minimal Costa surface is embedded in $\bfR^3$, but 
a corresponding CMC $1$ surface in $\bfH^3$ is not embedded.  
And the minimal trinoid is not embedded in $\bfR^3$, but a 
corresponding CMC $1$ surface in $\bfH^3$ is embedded.  
But again, this is a global question.  (Umehara and Yamada, Rosenberg, and 
others have 
been studying embeddedness of CMC $1$ surfaces in $\bfH^3$ \cite{chr1}, 
\cite{ET1}, \cite{LR}, \cite{uy1}.)  

The goals of this article are three-fold: 
\begin{enumerate}
\item to demonstrate some essential global 
differences between minimal surfaces in $\bfR^3$ and CMC $1$ surfaces 
in $\bfH^3$, for which graphics are helpful.  
\item to explain how period problems can be solved, without 
hiding the main ideas behind too many technical details.  Showing 
solvability invariably requires detailed arguments, but here we restrict 
ourselves to giving plausibility arguments, to make the main ideas 
more transparent.  
\item to show computer graphics of CMC $1$ surfaces in 
$\bfH^3$.   Such graphics are largely absent from an already significant 
amount of literature on the subject, but are useful for understanding 
the nature of these surfaces (in the same way that computer graphics have 
been useful in the study of minimal surfaces in $\bfR^3$).  Consider 
this example: the CMC $1$ Enneper cousin in $\bfH^3$ has an end at which the 
hyperbolic Gauss map (to be defined later) has an essential 
singularity, so it is difficult to imagine what this surface looks like, 
and no picture is in the literature.  Here we show its picture, and pictures 
of other fundamental examples as well.  
\end{enumerate}

In Section 1, we give an elementary description of CMC surfaces and the 
Lawson correspondence.  In Section 2, we describe hyperbolic space.  In 
Sections 3 and 4, we describe the Weierstrass representation for minimal 
surfaces in $\bfR^3$ and the Bryant representation for CMC $1$ surfaces 
in $\bfH^3$, respectively.  In Section 5, we describe the duals of CMC 
$1$ surfaces in $\bfH^3$, and then we describe a way to solve the period 
problems on dual surfaces in Sections 5 and 6.  
In the final section, we include an overview of known results about CMC 
$1$ surfaces in $\bfH^3$.  

\section{Defining CMC surfaces, and Lawson's correspondence}

Let $\Phi: \Sigma \to \bfR^3$ be an immersion of a $2$ dimensional surface 
$\Sigma$ into Euclidean 3-space $\bfR^3$.  The standard metric 
$dx_1^2+dx_2^2+dx_3^2$ (also written $\langle \cdot , \cdot \rangle$) on 
$\bfR^3$ induces a metric (the first fundamental form) $ds^2: 
T_p(\Sigma) \times T_p(\Sigma) \rightarrow \bfR$ on $\Sigma$, which 
is a bilinear map, 
where $T_p(\Sigma)$ is the tangent space at $p \in \Sigma$. If $(u,v)$ is 
a local coordinate of $\Sigma$, and if the basis 
$\{\Phi_u = (\frac{\partial \Phi}{\partial u})_p, \Phi_v = 
(\frac{\partial \Phi}{\partial v})_p \}$ is chosen for 
$T_p(\Sigma)$ (we can identify $T_p(\Sigma)$ with the plane in 
$\bfR^3$ tangent to $\Phi(\Sigma)$ at $\Phi(p)$), then the metric $ds^2$ is 
represented by the matrix 
\[ I_1 = \left( 
\begin{array}{cc}
E & F \\
F & G 
\end{array}
\right) = 
\left( 
\begin{array}{cc}
\langle \Phi_u, \Phi_u \rangle & \langle \Phi_u, \Phi_v \rangle \\
\langle \Phi_v, \Phi_u \rangle & \langle \Phi_v, \Phi_v \rangle 
\end{array}
\right) \; , \] that is, $ds^2(a\Phi_u+b\Phi_v,c\Phi_u+d\Phi_v)=
(a \; b)I_1(c \; d)^t$ for any $a,b,c,d \in \bfR$.  
Likewise, the second fundamental form of $\Phi$, a symmetric bilinear 
map from $T_p(\Sigma) \times T_p(\Sigma)$ to 
the normal space $\{r \cdot \vec{N}_p \; | \; r \in \bfR \} \approx \bfR$, 
where $\vec{N}$ is a local unit normal vector field to 
$\Phi(\Sigma)$ near $\Phi(p)$, is similarly represented by a matrix 
\[ I_2 = \left( 
\begin{array}{cc}
l & m \\
m & n 
\end{array}
\right) = - 
\left( 
\begin{array}{cc}
\langle \vec{N}_u, \Phi_u \rangle & \langle \vec{N}_u, \Phi_v \rangle \\
\langle \vec{N}_v, \Phi_u \rangle & \langle \vec{N}_v, \Phi_v \rangle 
\end{array} \right) \; . \]  

The linear map $(a \; b)^t \to I_1^{-1} I_2 (a \; b)^t$ represents the 
shape operator taking $\nu = a \Phi_u + b \Phi_v \in T_p(\Sigma)$ to 
$-D_\nu N \in T_p(\Sigma)$, where $D$ is the directional derivative of 
$\bfR^3$.  The eigenvalues $k_1, k_2$ and corresponding eigenvectors of 
the shape operator are the principal curvatures and principal 
curvature directions of the surface $\Phi(\Sigma)$ at 
$\Phi(p)$.  Although a bit imprecise, intuitively the principal curvatures 
tell us the maximum and minimum amounts amounts of bending 
of the surface toward the normal $\vec{N}$ at $\Phi(p)$, 
and the principal curvature directions tell us the directions in 
$T_p(\Sigma)$ of those maximal and minimal bendings.  

\begin{definition}
The determinant and the trace of the shape operator are 
the {\em Gaussian curvature} $K$ and the {\em 
mean curvature} $H$, respectively.  Thus \[ K = k_1k_2 = 
\frac{ln-m^2}{EG-F^2} \; , \; \; \; \; \; \mbox{and} \; \; \; \; \;  
H = \frac{k_1 + k_2}{2} = \frac{Gl-2Fm+En}{2(EG-F^2)} \; . \]  
A surface is {\em CMC} if $H$ is constant, 
and is {\em minimal} if $H$ is identically zero.  
\end{definition}

It is a fact that one can locally always choose $(u,v)$ so that $z=
u+iv$ becomes a conformal coordinate on $\Sigma$ with respect to $ds^2$ (i.e. 
$F \equiv 0$, $E \equiv G$).  
Hence from now on we assume the coordinates are chosen 
this way and furthermore that $\Sigma$ is a Riemann surface.  

The mean curvature $H$ at $\Phi(p)$ can be equivalently defined 
as the average of the normal curvatures $-\langle \nu, 
D_\nu \vec{N} \rangle$ in all 
tangent directions $\nu \in T_p(\Sigma)$.  (Intuitively, the normal curvature 
measures the rate at which the surface bends toward $\vec{N}$, in the 
direction $\nu$.)  Thus a minimal surface has 
average normal curvature zero at every point, and this suggests a physical 
interpretation: 
\begin{quote}
\cite{HoMe2}: 
``Loosely speaking, one imagines the surface as made up of very 
many rubber bands, stretched out in all directions; on a minimal surface the 
forces due to the rubber bands balance out, and the surface does not need to 
move to reduce tension.''  
\end{quote}
To say this more rigorously, suppose $U$ is a compact domain in $\Sigma$, and 
define a {\em smooth boundary-fixing 
variation} of the immersion $\Phi (U)$ to be a $C^\infty$ 
map $\Phi_t: (-1,1) \times U \to \bfR^3$ with three properties: 
\begin{enumerate}
\item $\Phi_t(\cdot): U \to \bfR^3$ is an immersion for all 
$t \in (-1,1)$, 
\item $\Phi_0 = \Phi$ on $U$, 
\item $\Phi_t|_{\partial U} = \Phi|_{\partial U}$ for all 
$t \in (-1,1)$.  
\end{enumerate} 
Note that 
$\mbox{Area}(\Phi_t(U)) = \int_\Sigma dA_t$, where $dA_t$ is the volume 
element (the area $2$ form) 
of the metric induced by the immersion $\Phi_t$.  It turns out that the 
first variation formula for smooth boundary-fixing variations is then 
\[ \left. \frac{d}{dt} \mbox{Area}(\Phi_t(U)) \right|_{t=0} = 
- \int_\Sigma \left\langle H \vec{N}, 
\left. (\Phi_t)_* \frac{\partial}{\partial 
t} \right|_{t=0} \right\rangle dA_{0} \; . \]  
In particular, minimal surfaces (with $H \equiv 0$) are area critical for all 
compact domains $U$ -- and we could 
have defined them this way.  (Actually, when $U$ is small enough, 
not only is $\Phi(U)$ critical for area, it is also the 
unique least-area surface with boundary $\partial \Phi(U)$, hence 
a "minimal" surface.)  

Similarly, a nonminimal CMC surface could be defined 
as an immersion $\Phi$ such that for every compact 
domain $U$, $\Phi(U)$ is critical for area amongst all 
smooth boundary-fixing variations that keep the volume on one side of the 
surface unchanged.  (The derivative of volume is 
$\int_U \langle \vec{N}, \left. (\Phi_t)_* \frac{\partial}{\partial 
t} \right|_{t=0} \rangle dA_t$, so if the volume is unchanging and so 
$\int_U \langle \vec{N}, \left. (\Phi_t)_* \frac{\partial}{\partial 
t} \right|_{t=0} \rangle dA_t=0$, and if $H$ is constant, then 
$\left. \frac{d}{dt} \mbox{Area}(\Phi_t(U)) \right|_{t=0} =0$ \cite{BCE}.)  

This is why minimal and CMC surfaces model 
physical soap films, which always move to minimize area.  
Minimal surfaces model soap films 
not enclosing bounded pockets of air, as such films 
are area minimizing for all boundary-fixing 
variations.  Nonminimal CMC surfaces model soap films enclosing 
bounded pockets of air, as such films are area minimizing only 
for variations that keep the air pockets' volumes fixed.  

To define CMC surfaces in $\bfH^3$, we can proceed in the same way as 
the Euclidean case.  We only need to replace 
$\Phi_u$ and $\Phi_v$ by linear differentials $(\frac{\partial}
{\partial u})_p$ and $(\frac{\partial}{\partial v})_p$, and 
$\vec{N}_u$ and 
$\vec{N}_v$ by $\nabla_{\frac{\partial}{\partial u}} \vec{N}$ and 
$\nabla_{\frac{\partial}{\partial u}} \vec{N}$, 
where $\nabla$ is the Riemannian connection for $\bfH^3$.  The shape 
operator becomes $\nu \to -\nabla_{\nu} N$.  All the 
variational properties in the Euclidean case also hold 
when the ambient space is $\bfH^3$.  And we can choose 
conformal coordinates $z=u+iv$ in the $\bfH^3$ case as well, so we again 
assume $\Sigma$ is a Riemann surface.  

\begin{figure}
\centerline{
        \hbox{
		\psfig{figure=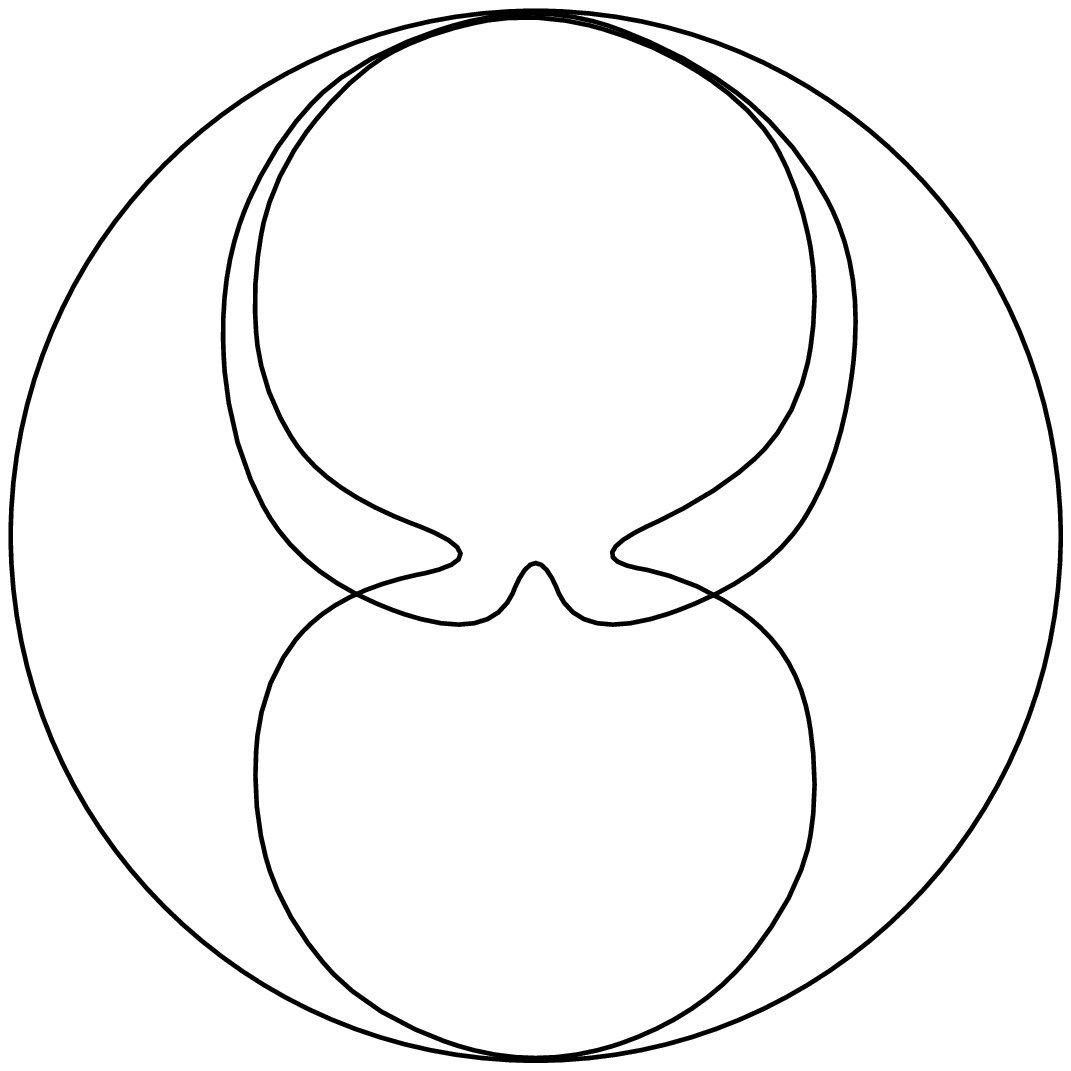,width=2in}
		\psfig{figure=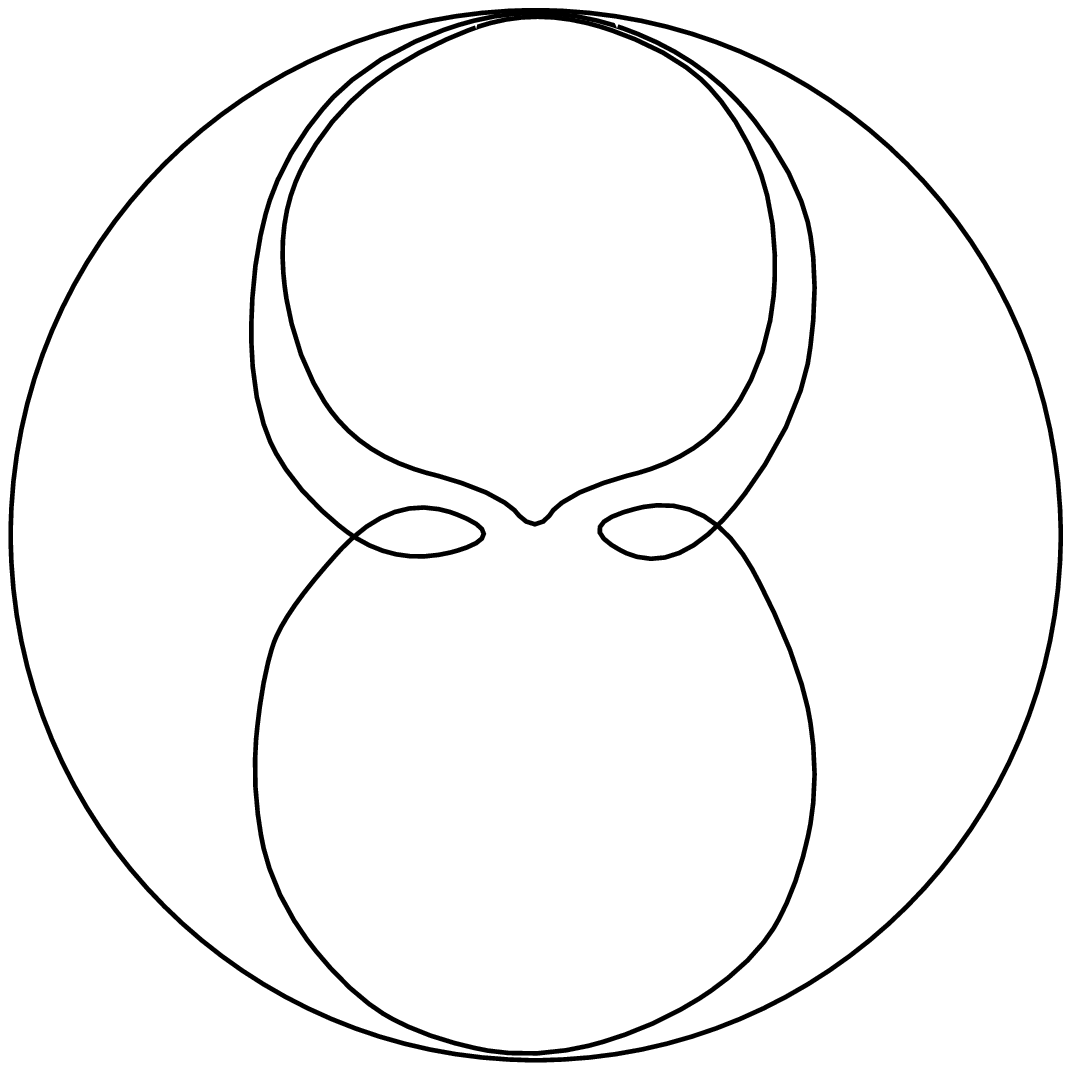,width=2in}
		\psfig{figure=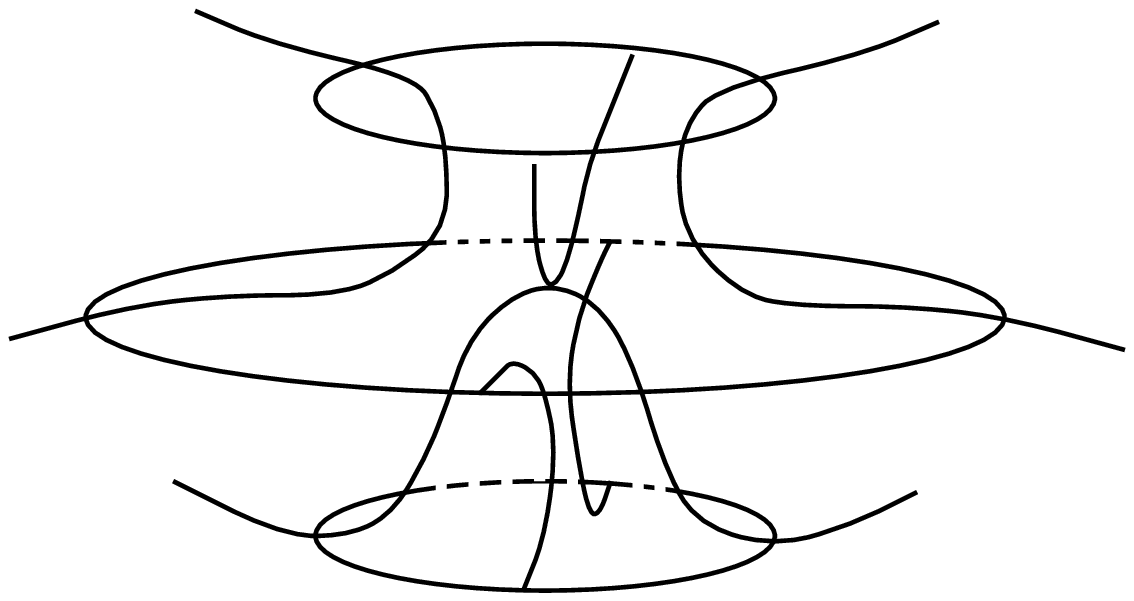,width=2in}
	}
}
\caption{CMC $1$ Costa cousin dual in $\bfH^3$, proven to exist by 
Costa and Sousa Neto when they established that the minimal Costa 
surface in $\bfR^3$ is symmetric and nondegenerate \cite{CN}.  Rather 
than showing graphics of this surface, we show two vertical cross sections by 
which the surface is reflectionally symmetric (including the "circles" at 
infinity), and a schematic of the central portion of the surface.}
\label{figure14.5}
\end{figure}

\begin{figure}
\vspace{0.5in}
\centerline{
        \hbox{
		\psfig{figure=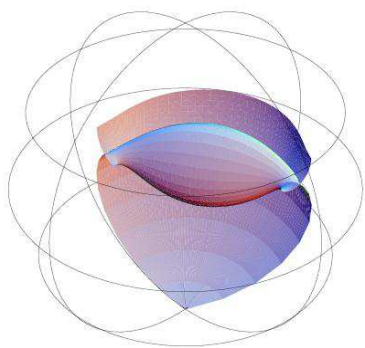,width=2in}
	}
}
\caption{
$5$ ended CMC $1$ surface in $\bfH^3$, found by Umehara and 
Yamada in \cite{uy1}. This surface has the same symmetries and 
end behavior as the minimal surface on the left hand side of Figure 11.  
Here we show only one of six congruent disks that form the surface.  
The full surface is constructed by reflections across planes 
containing boundary curves of the disk shown here.}
\label{figure17}
\end{figure}

Regarding the Lawson correspondence between
minimal surfaces in $\bfR^3$ and CMC $1$ surfaces in $\bfH^3$, the essential 
ingredient is the fundamental theorem 
of surface theory, telling us that on a simply connected region $U$ of 
$\Sigma$ there exists an 
immersion with fundamental forms $I_1$ and $I_2$ if and only if 
$I_1$ and $I_2$ satisfy the Gauss and Codazzi equations.  (This is true 
regardless of whether the ambient space is $\bfR^3$ or $\bfH^3$, but the 
Gauss and Codazzi equations are not the same in the two cases.)  So if 
$\Phi : U \to \bfR^3$ is a conformal minimal immersion, then its 
fundamental forms $I_1$ and $I_2$ satisfy the 
Gauss and Codazzi equations for surfaces in $\bfR^3$, and one can then check 
that $\tilde{I}_1=I_1$ and $\tilde{I}_2=I_1+I_2$ satisfy the Gauss and 
Codazzi equations for surfaces in $\bfH^3$, so there exists an immersion 
$\tilde{\Phi} : U \rightarrow
\bfH^3$ with fundamental forms $\tilde{I}_1$ and 
$\tilde{I}_2$.  Since $\mbox{tr}(I_1^{-1}I_2)=0$, we have 
$\mbox{tr}(\tilde{I}_1^{-1} \tilde{I}_2)=1$, so 
$\tilde{\Phi}$ is CMC $1$.  And since $I_1 = \tilde{I}_1$, $\tilde{\Phi}$ is 
also conformal, and $\Phi(U)$ and $\tilde{\Phi}(U)$ are 
isometric.  This is Lawson's correspondence.  

\begin{quote}
{\small 

We digress here to explain the Lawson correspondence in more detail: Let 
$M^3(\bar{K})$ be the $3$ dimensional space form with constant section 
curvature $\bar{K}$ (e.g. $M^3(0)=\bfR^3$, $M^3(-1)=\bfH^3$, 
$M^3(1)=S^3$).  For an immersion $\Phi:U \to M^3(\bar{K})$ with induced 
metric $\langle \cdot , \cdot \rangle$ and connection $\nabla$ and Gaussian 
curvature $K$ and shape operator $A$, the Gauss and Codazzi equations are 
satisfied:  
\[ K - \bar{K} = \det(A) \; , \; \; \; 
\langle A([X,Y]),Z \rangle = 
\langle \nabla_X A(Y),Z \rangle - 
\langle \nabla_Y A(X),Z \rangle \]  for all smooth vector fields 
$X$, $Y$, and $Z$ in the tangent space of $U$.  
Assume $\Phi$ is CMC $H$, so $H = \mbox{tr}(A)$ is constant.  
Now choose any $c \in \bfR$ and define 
\[ \tilde{A} = A+c \cdot (id.) \; , \; \; \tilde{K} = 
\bar{K} - 2 c \mbox{tr}(A) - c^2 \; . \] 
Then the Codazzi equation clearly still holds when $A$ is replaced by 
$\tilde{A}$:  \[ \langle \tilde{A}([X,Y]),Z \rangle = 
\langle A([X,Y]),Z \rangle + c \langle [X,Y],Z \rangle = \langle 
\nabla_X A(Y),Z \rangle - \]\[ 
\langle \nabla_Y A(X),Z \rangle + c \langle \nabla_X Y,Z \rangle - 
c \langle \nabla_Y X,Z \rangle = \langle \nabla_X 
\tilde{A}(Y),Z \rangle - \langle \nabla_Y \tilde{A}(X),Z \rangle \; . \] 
The Gauss equation in $M^3(\tilde{K})$ also holds with $\tilde{A}$ 
(note that $K$ is intrinsic and does not change): 
\[ K-\tilde{K} = K - (\bar{K} -2c\mbox{tr}(A) - c^2) = \det(A) + 
2c\mbox{tr}(A) + c^2 \]\[ = \det(A+c(id.)) = \det(\tilde{A}) \; . \] 
Therefore there exists 
an immersion $\tilde{\Phi}: U \to M^3(\tilde{K})$ with metric 
$\langle \cdot , \cdot \rangle$ and shape operator $\tilde{A}$, and 
$\tilde{\Phi}(U)$ is isometric to $\Phi(U)$.  
As the mean curvature of $\tilde{\Phi}(U)$ is 
\[ \tilde{H} = \mbox{tr}(\tilde{A}) = 
\mbox{tr}(A) + c = H + c \; , \] this 
demonstrates the Lawson correspondence between a CMC $H$ surface in 
$M^3(\bar{K})$ and a CMC $(H+c)$ surface in $M^3(\bar{K} - 2 c H - c^2)$.  
In particular, when $H=\bar{K}=0$ and $c=1$, we have the 
correspondence between minimal surfaces in $\bfR^3$ and CMC $1$ surfaces in 
$\bfH^3$.  

}
\end{quote}

\section{Hyperbolic space}

\begin{figure}
\vspace{0.5in}
\centerline{
        \hbox{
		\psfig{figure=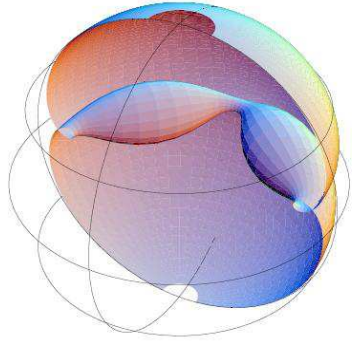,width=2in}
	}
}
\caption{CMC $1$ genus $0$ surface in $\bfH^3$ with two embedded ends, 
which is not a surface of revolution.  By Schoen's result, such a surface 
cannot exist as a minimal surface in $\bfR^3$, so this surface is not 
the cousin nor the dual cousin of any minimal surface with finite 
total curvature in $\bfR^3$.  Its 
existence is proven in \cite{uy1}.  One of two congruent pieces of the 
surface is shown here.  This example, and the examples in 
Figures 7 and 13, show that 
the converse of the \cite{ruy1} result does not hold.}
\label{figure15}
\end{figure}

\begin{figure}
\vspace{0.5in}
\centerline{
        \hbox{
		\psfig{figure=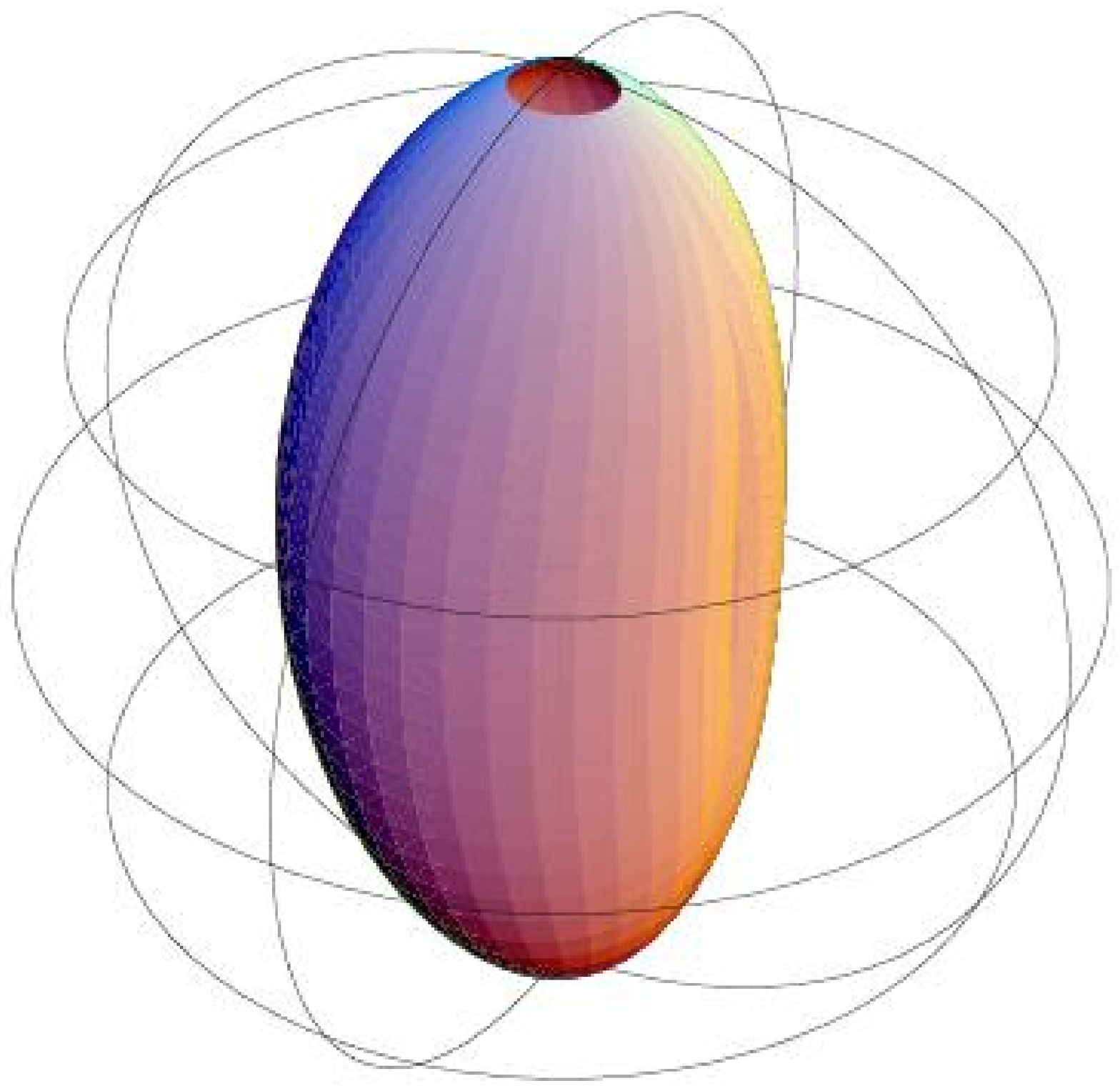,width=2in}
		\hspace{0.5in}
		\psfig{figure=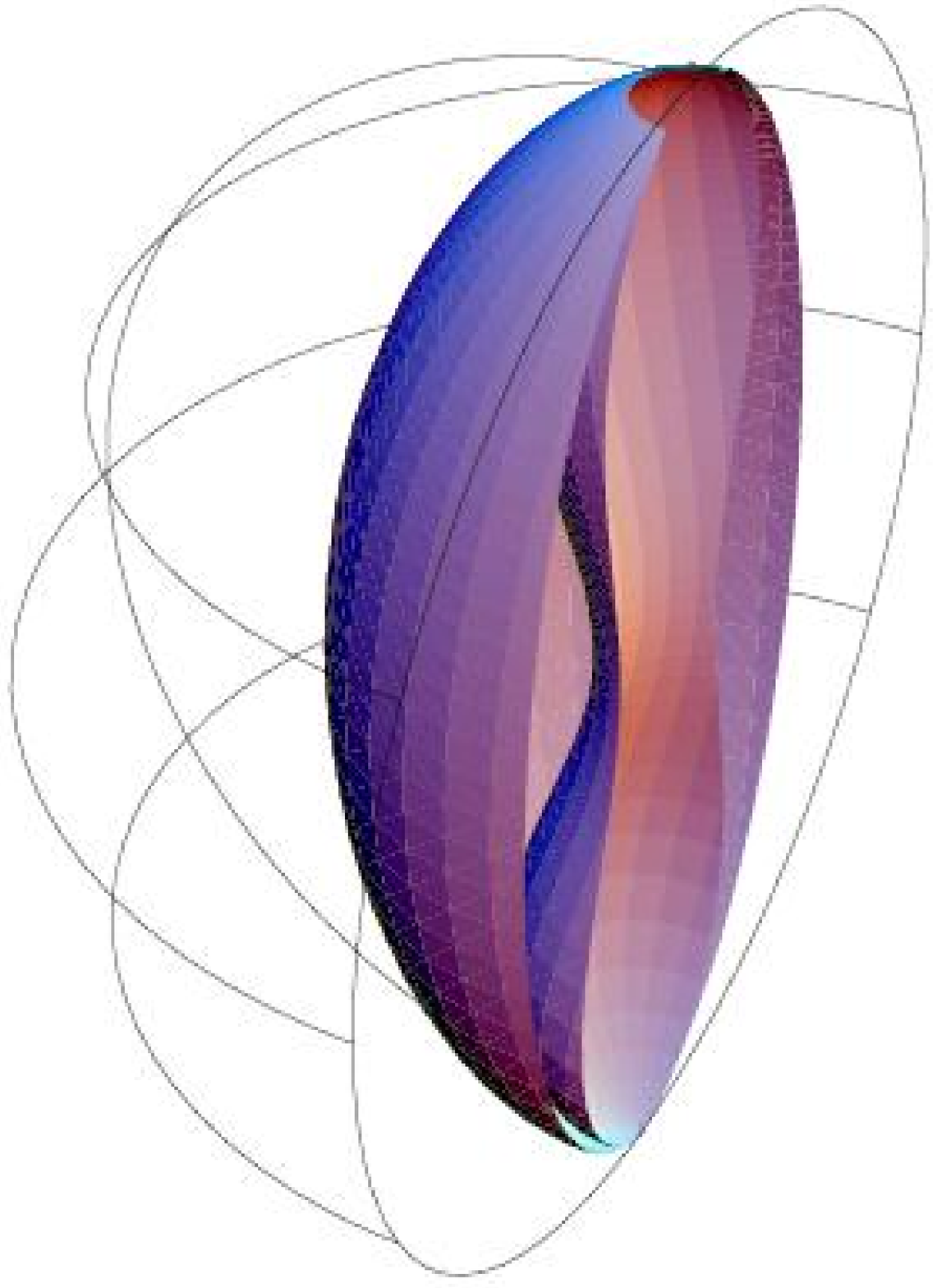,width=1.5in}}}
		\vspace{0.5in}
		\centerline{\hbox{
		\psfig{figure=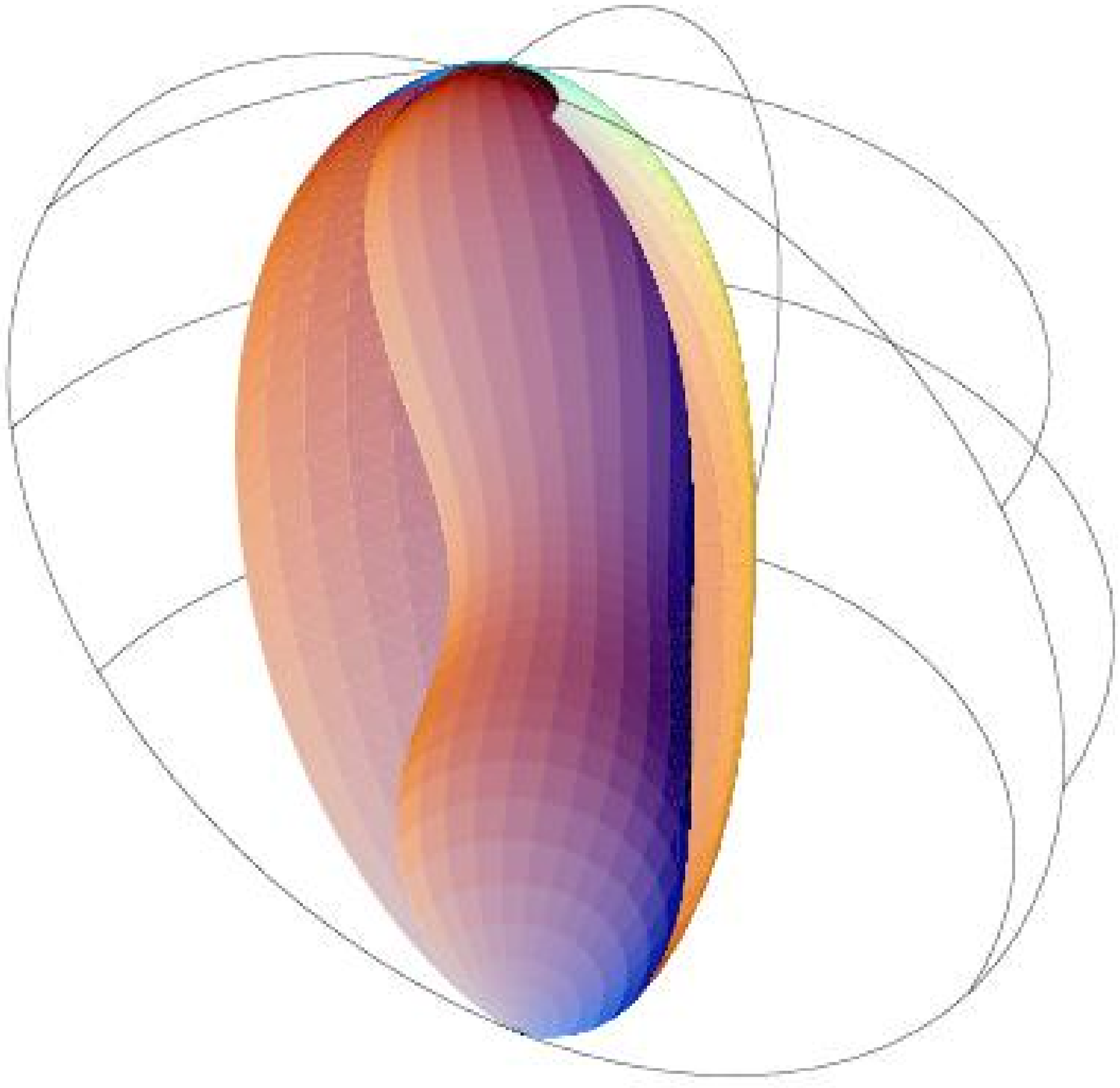,width=1.7in}
		\hspace{0.5in}
		\psfig{figure=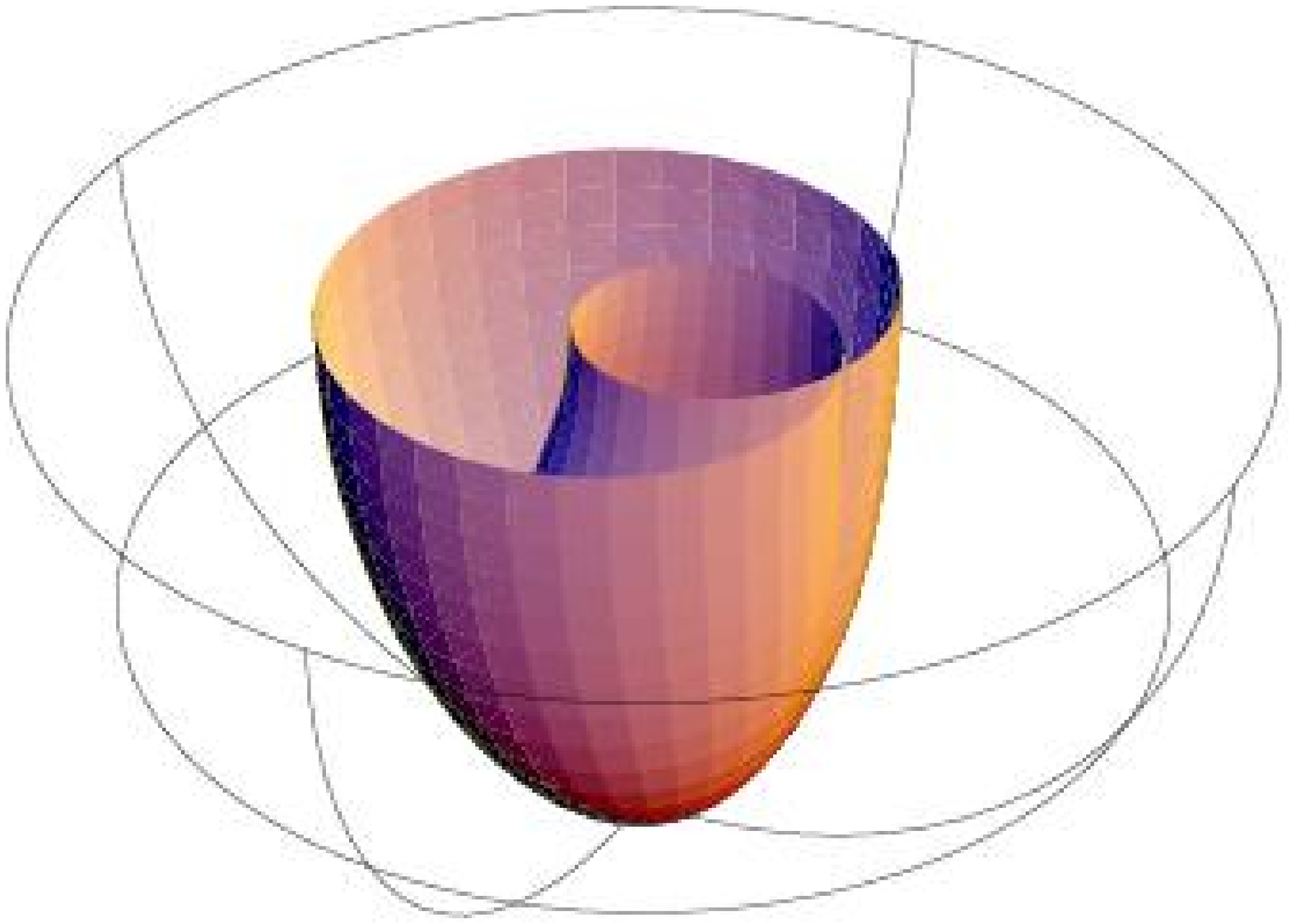,width=2in}
	}
}
\caption{CMC $1$ genus $0$ surface in $\bfH^3$ with two nonembedded 
ends (each of winding order two), proven to exist in \cite{uy1}.  
This surface again (like the surface in Figure 6) 
has no corresponding minimal surface in $\bfR^3$.  
Three partial views of the surface are also shown here.  }
\label{figure16}
\end{figure}

Hyperbolic $3$-space $\bfH^3$ is the unique simply-connected 
$3$ dimensional complete Riemannian manifold with 
constant sectional curvature $-1$, but it can be described by a 
variety of models, each with their own advantages.  Here we describe 
the two models we will need: the Poincare ball model, 
convenient for showing computer graphics; and the Hermitian 
matrix model, the one Bryant chose for his representation.  

The Poincare model is the $3$ dimensional Euclidean unit ball 
\[ B^3 = \{(x_1,x_2,x_3) \in \bfR^3 \; | \; \; 
x_1^2+x_2^2+x_3^2 < 1\} \] 
with the metric \[ ds^2 = (\frac{2}{1 - x_1^2-x_2^2-x_3^2})^2 
(dx_1^2+dx_2^2+dx_3^2) . \]  This metric is
conformal to the Euclidean metric and so 
angles are the same as Euclidean angles, that is, angles are just as our 
"Euclidean" eyes tell us they are.
However, distances are clearly not Euclidean.  With this metric 
the unit ball is complete, simply-connected, and has 
constant sectional curvature $-1$, so it is $\bfH^3$.  

The geodesics in the Poincare model are segments of Euclidean
lines and circles that intersect the boundary "sphere at infinity" $\partial 
B^3$ at right angles.  The CMC hyperbolic planes with $H = 0$ (resp. 
hyperspheres with $|H| < 1$, spheres with $|H| > 1$, horospheres with $|H|=1$) 
are the intersections of $B^3$ with Euclidean spheres and planes in $\bfR^3$ 
that meet $\partial B^3$ orthogonally (resp. have nonempty nontangential 
intersection with $\partial B^3$, are in the interior of $B^3$, are 
in $\overline{B^3}$ and are tangent to $\partial B^3$ at one point).  

Now we describe the Hermitian model, which we will see is convenient 
for describing the isometric motions of $\bfH^3$.  
We first recall that the $6$ dimensional 
Lie group $SL(2,\bfC)$ is the set of all $2 \times 2$ matrices with complex
entries and determinant 1, with matrix multiplication as the group 
operation.  We also mention, as we use it later, that 
$SU(2)$ is the $3$ dimensional subgroup of matrices 
$\mu \in SL(2,\bfC)$ such that $\mu \mu^*$ is the identity, where 
$\mu^* = \bar{\mu}^t$, or equivalently, 
\[ \mu = \left(
\begin{array}{cc}
p & -\bar{q} \\ q & \bar{p}
\end{array} \right) \; , \]
for some $p$, $q \in \bfC$ with $|p|^2 + |q|^2 = 1$.

The set 
\[ \{A A^* \; | \; A \in SL(2,\bfC) \} \;  \]
will form the Hermitian model for $\bfH^3$, which consists of 
the Hermitian symmetric matrices with determinant $1$, that is, 
matrices of the form 
\[ \left( \begin{array}{cc}
a_{11} & a_{12} \\ \overline{a_{12}} & a_{22} 
\end{array} \right) \; , \] where $a_{12} \in \bfC$ and $a_{11}, 
a_{22} \in \bfR$ and $a_{11} a_{22} - a_{12} \overline{a_{12}} = 1$.  
Such matrices can be bijectively mapped to points 
\[ \left( \frac{a_{12}+\overline{a_{12}}}{2+a_{11}+a_{22}}, 
\frac{i(\overline{a_{12}}-a_{12})}{2+a_{11}+a_{22}}, 
\frac{a_{11}-a_{22}}{2+a_{11}+a_{22}} \right) \] in the Poincare 
model, so if this set of matrices is given the metric so that 
this map is an isometry, then it becomes $\bfH^3$.  
We can then use the group $SL(2,\bfC)$ to describe 
the isometries of $\bfH^3$.  A matrix $h \in SL(2,\bfC)$ acts 
isometrically on $\bfH^3$ via, for each matrix $x$ in the Hermitian model, 
\[x \to h \, x \, h^* \; . \] 

\section{Weierstrass Representation}

In the Weierstrass representation for complete oriented 
minimal surfaces of finite total curvature (i.e. $\int_\Sigma -K dA 
< + \infty$) in $\bfR^3$, we incorporate that these surfaces are 
conformally equivalent to compact Riemann surfaces 
with finitely many points removed \cite{O}:  

\begin{figure}
\vspace{0.5in}
\centerline{
        \hbox{
		\psfig{figure=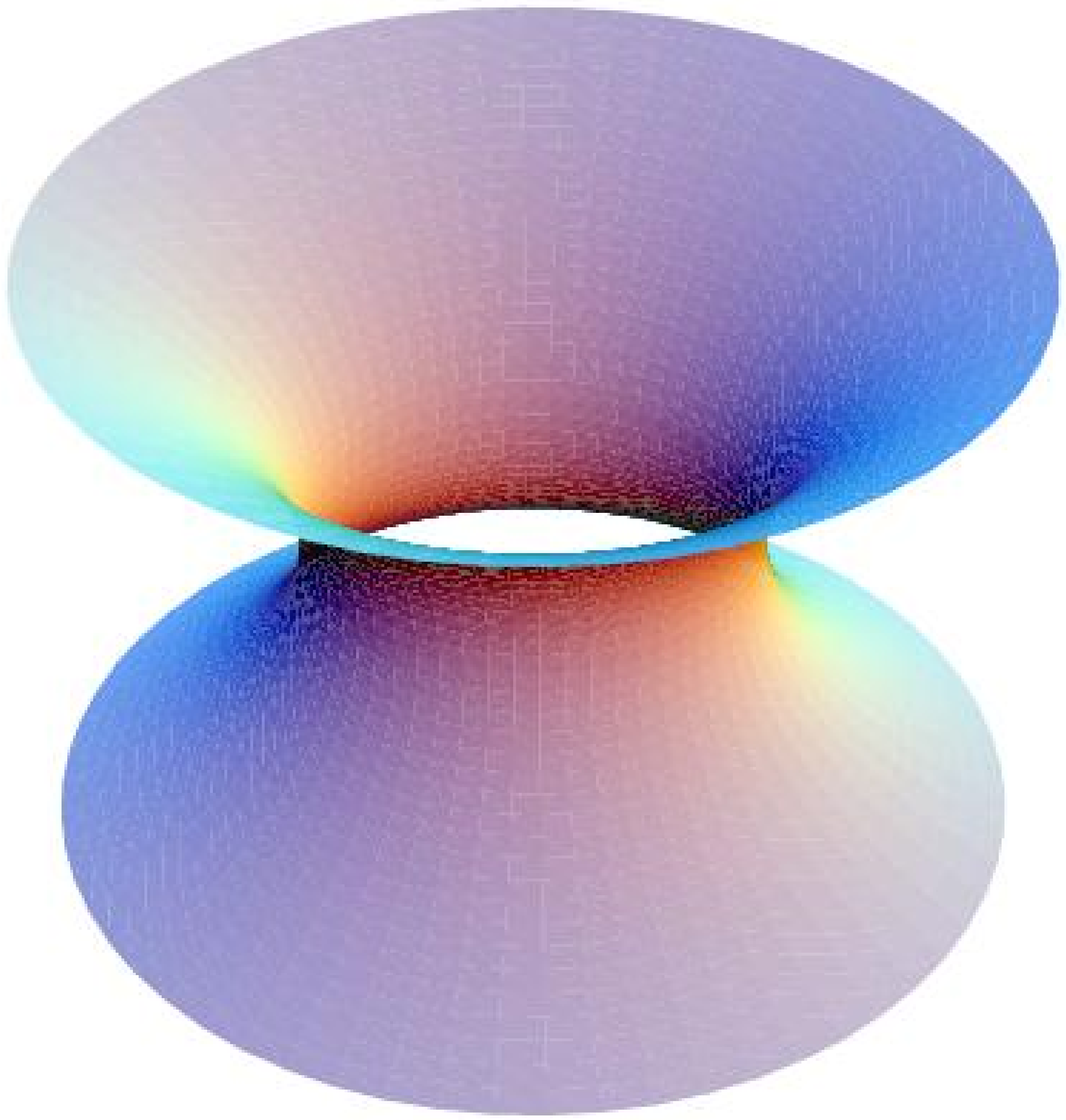,width=2in}
		\hspace{0.5in}
		\psfig{figure=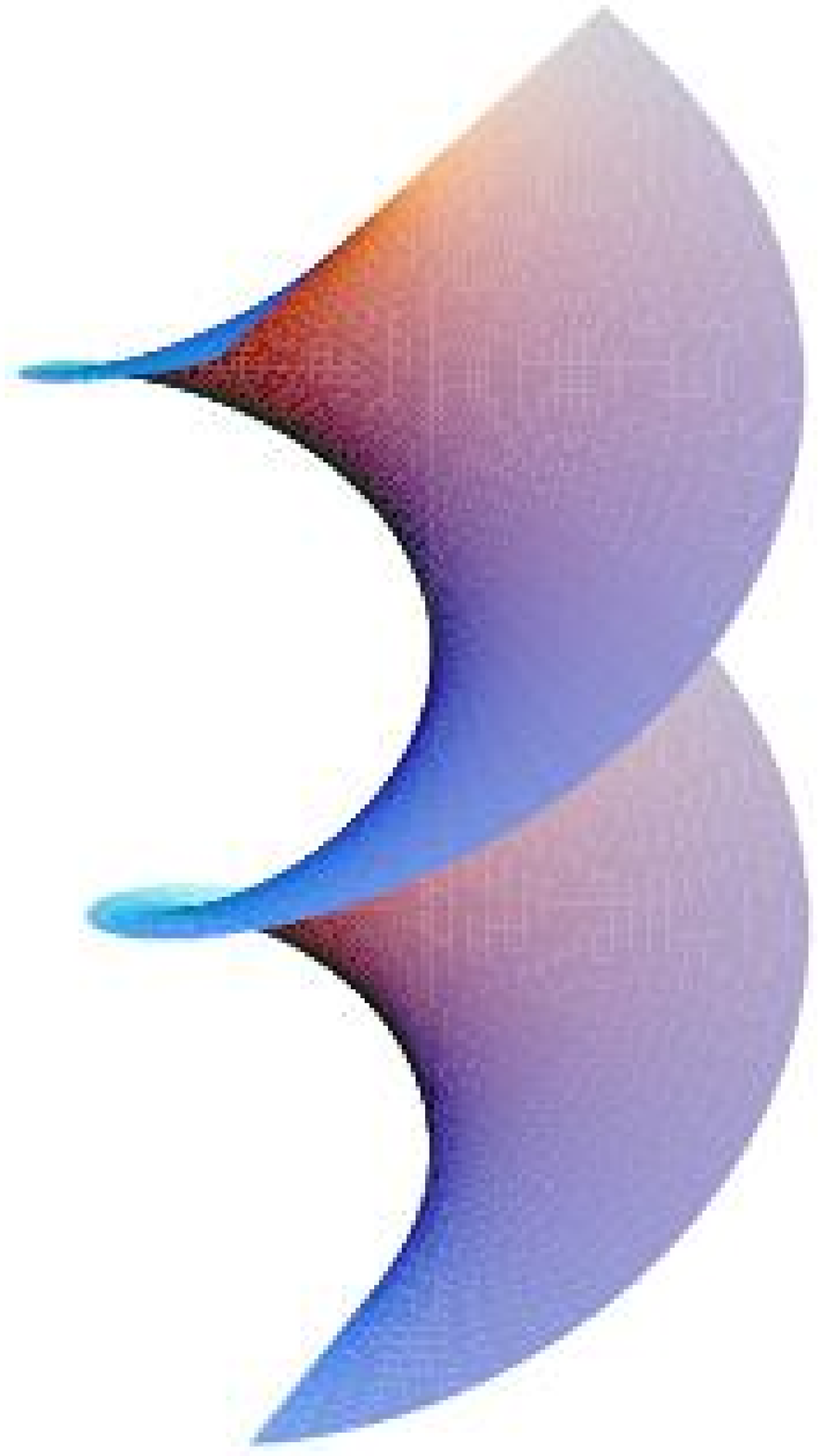,width=1.5in}
	}
}
\caption{Minimal catenoid and helicoid in $\bfR^3$.  (Note that the ends 
have been cut away, and will also be cut away in all subsequent figures.)  
Schoen \cite{S} has shown that any complete connected 
finite total curvature minimal immersion (not necessarily embedded) with two 
embedded ends in $\bfR^3$ must be a catenoid.}
\label{figure0}
\end{figure}

\begin{figure}
\vspace{0.5in}
\centerline{
        \hbox{
		\psfig{figure=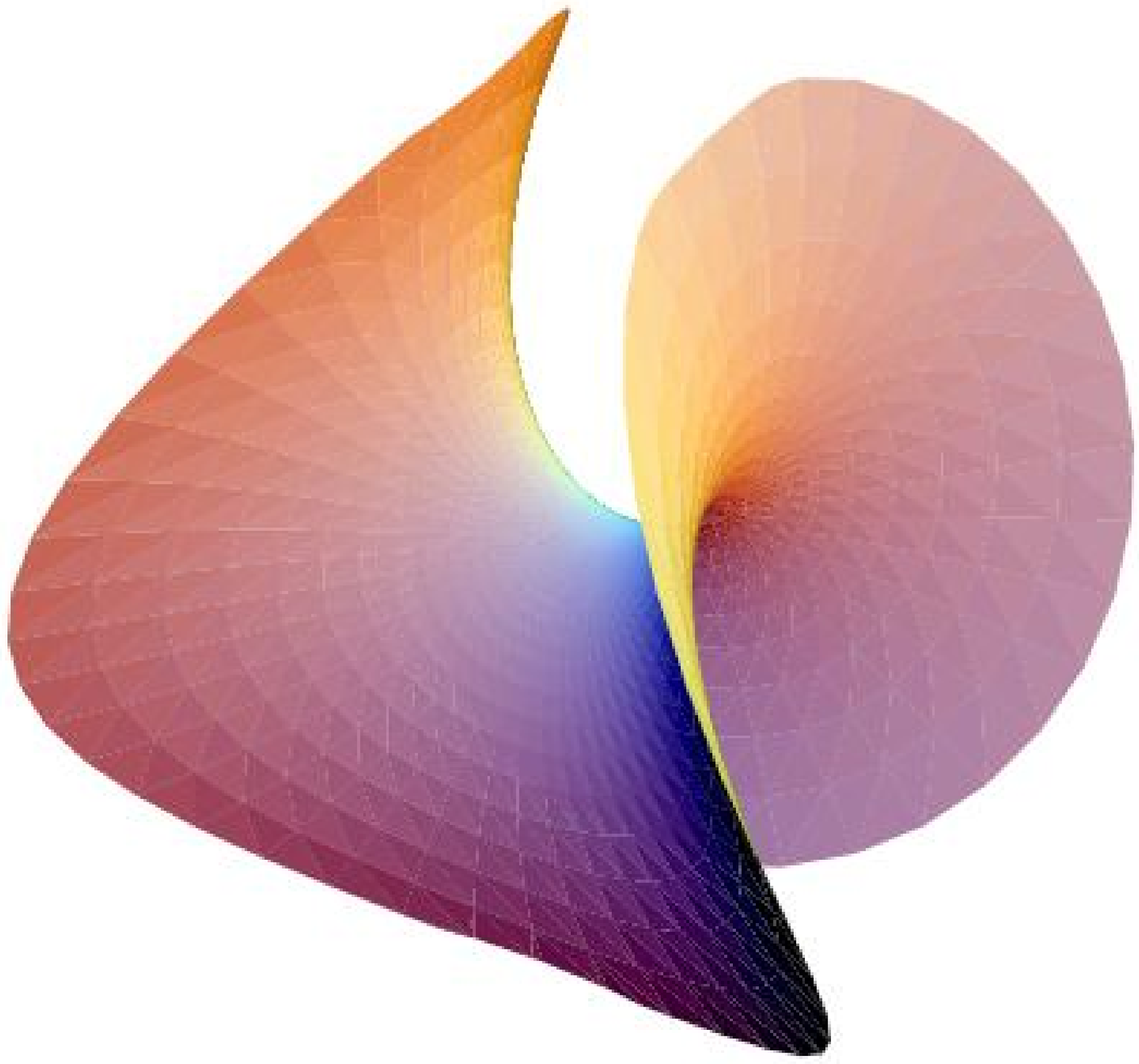,width=2in}
		\hspace{0.5in}
		\psfig{figure=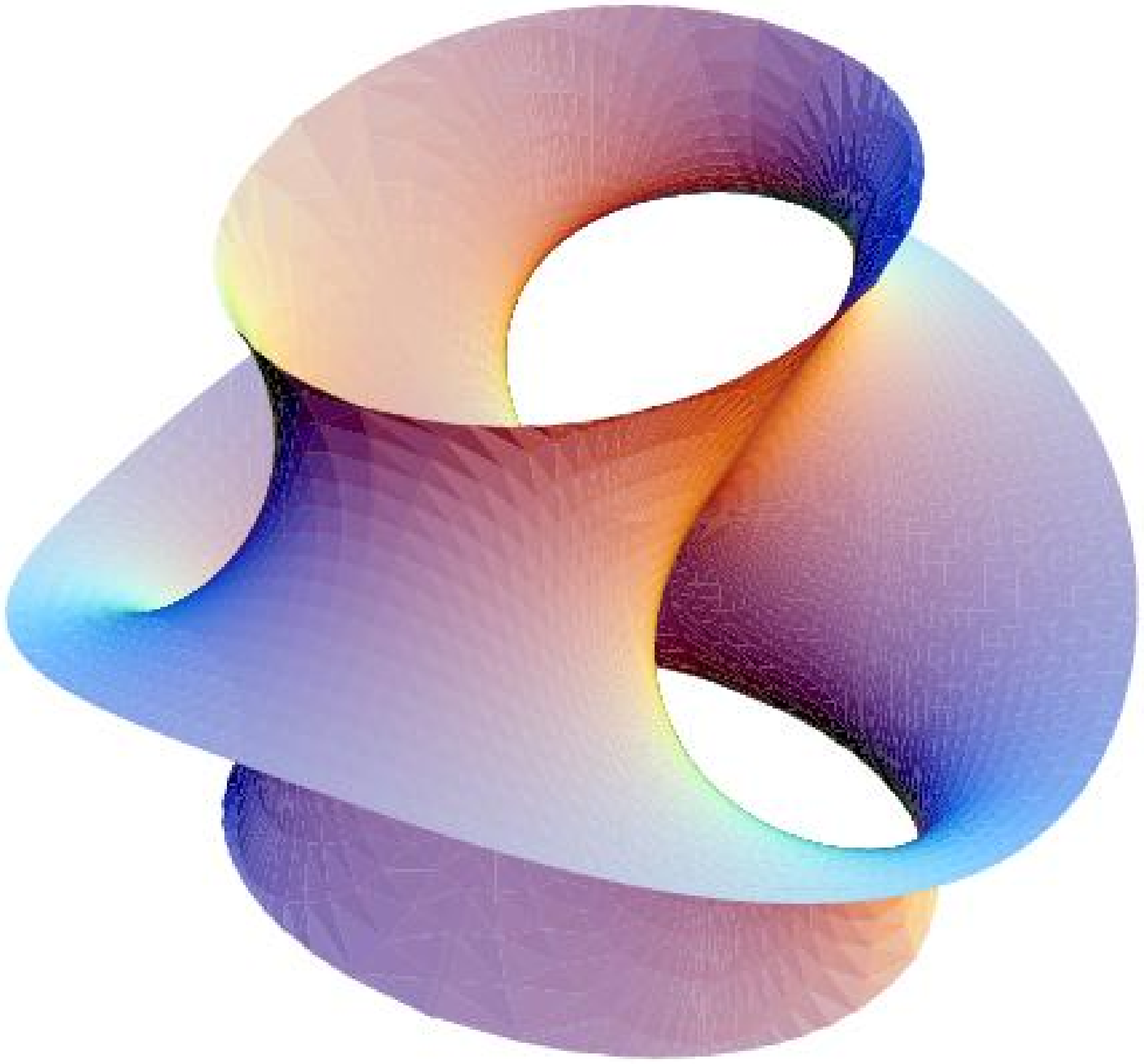,width=2in}
	}
}
\caption{Minimal Enneper surface and Costa surface in $\bfR^3$, the 
Costa surface was proven to exist by Costa and proven to 
be embedded in \cite{HoMe}.}
\label{figure1}
\end{figure}

\begin{figure}
\centerline{
        \hbox{
		\psfig{figure=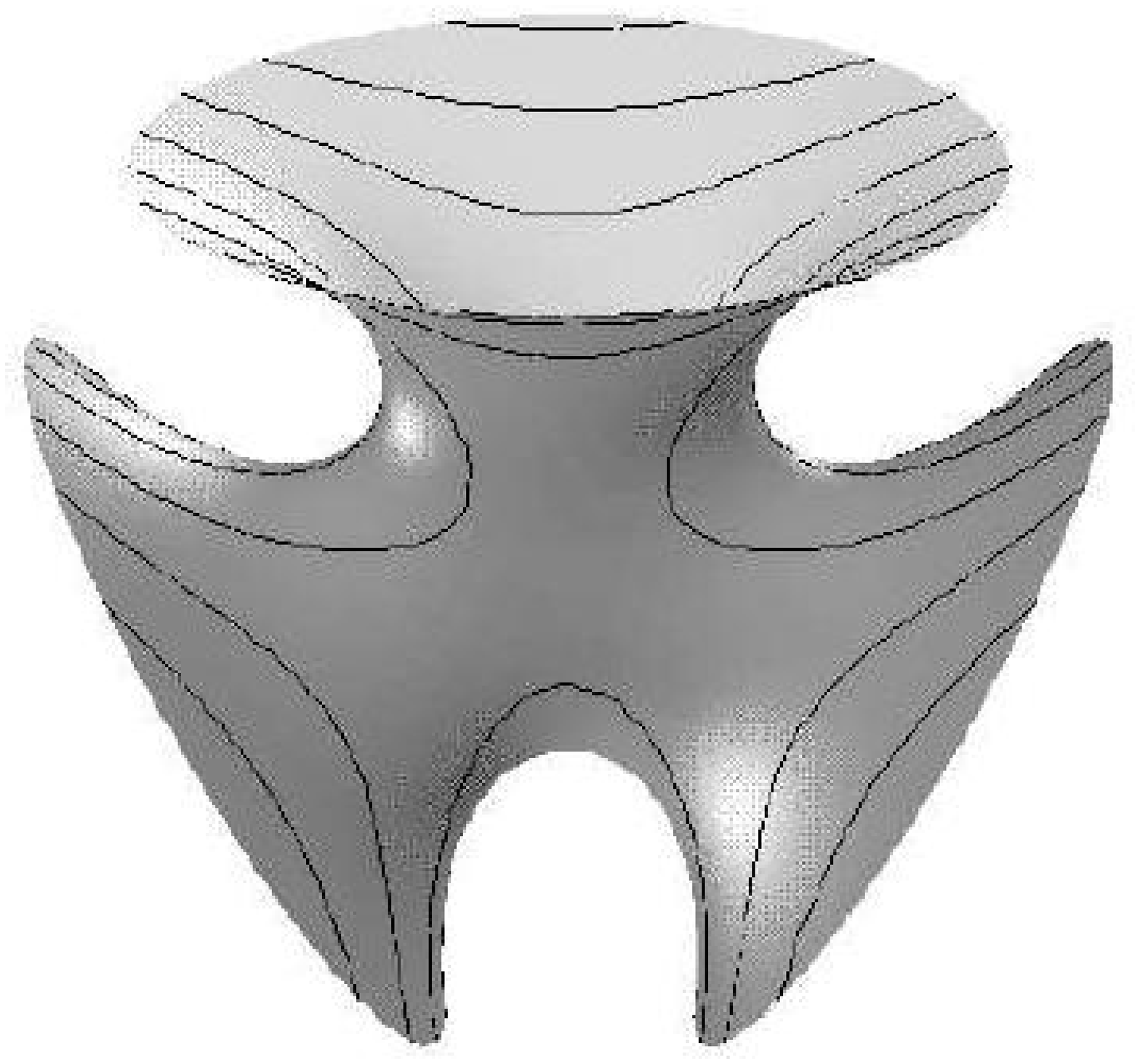,width=2.2in}
		\psfig{figure=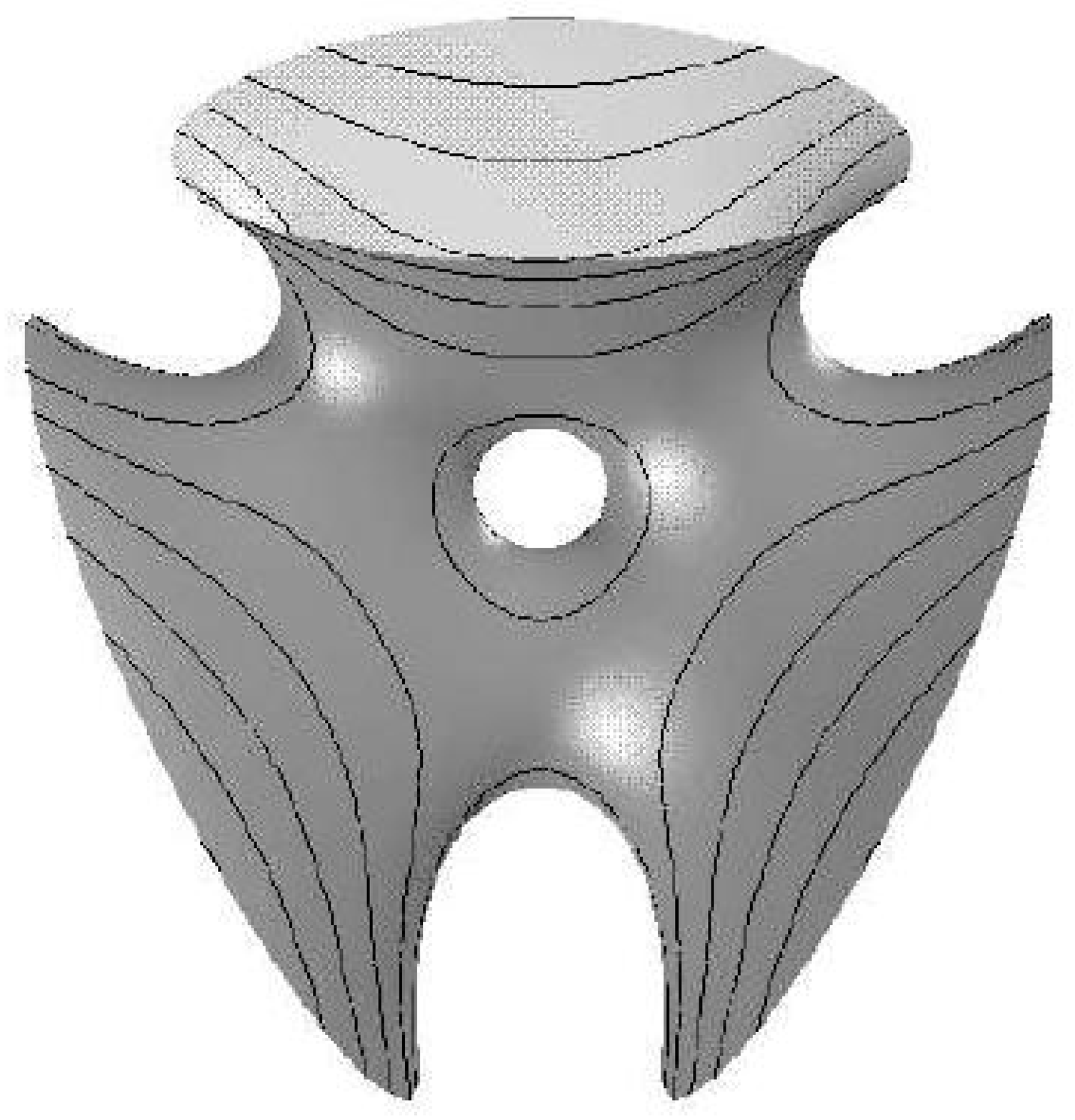,width=2.2in}
	}
}
\caption{Minimal trinoid (proven to exist in \cite{JM}) and genus 
$1$ trinoid (proven to exist in \cite{BR}) in $\bfR^3$.  These graphics 
were made using \cite{MESH} software.}
\label{figure4}
\end{figure}

\begin{figure}
\centerline{
        \hbox{
		\psfig{figure=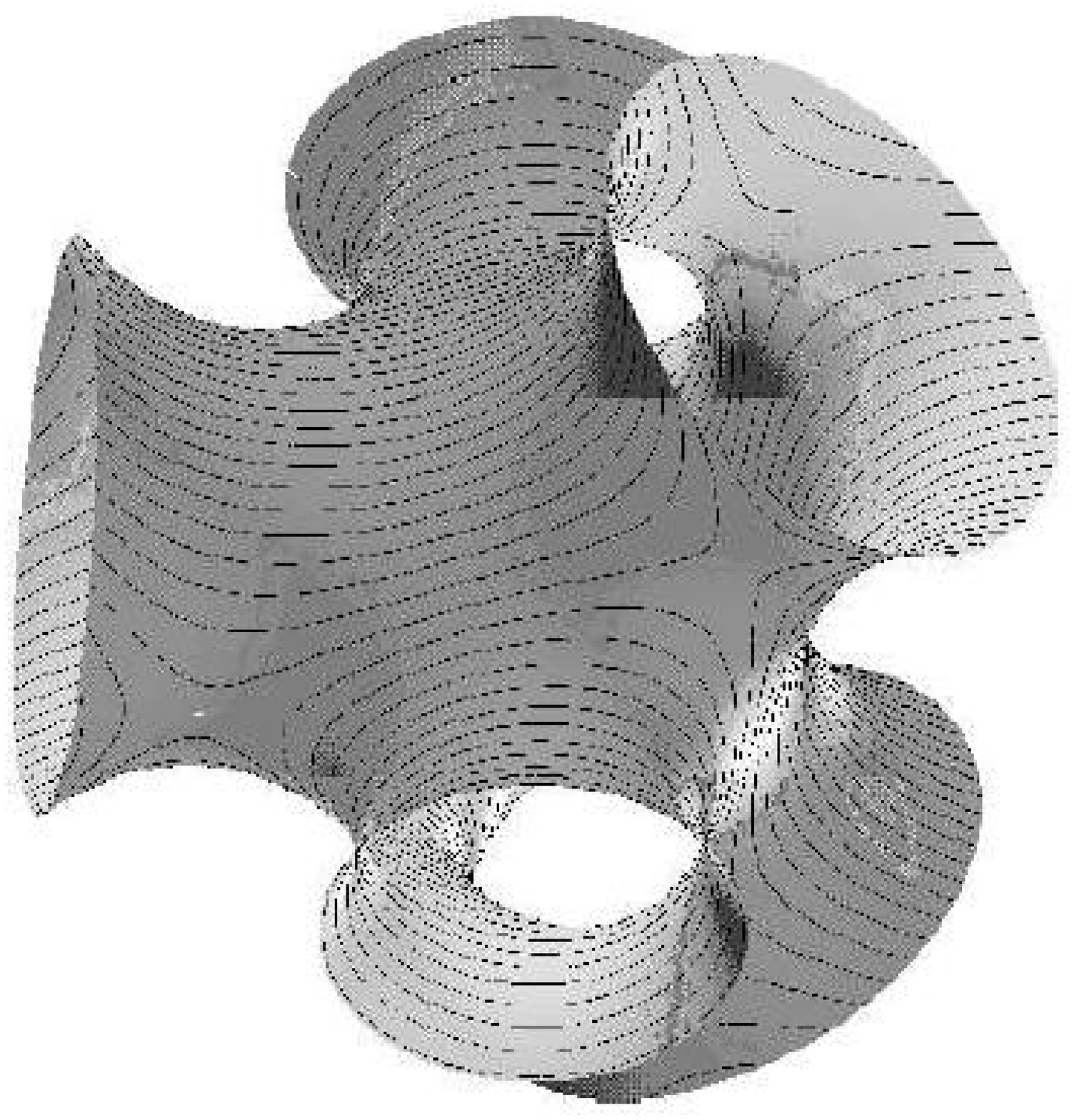,width=2.2in}
		\psfig{figure=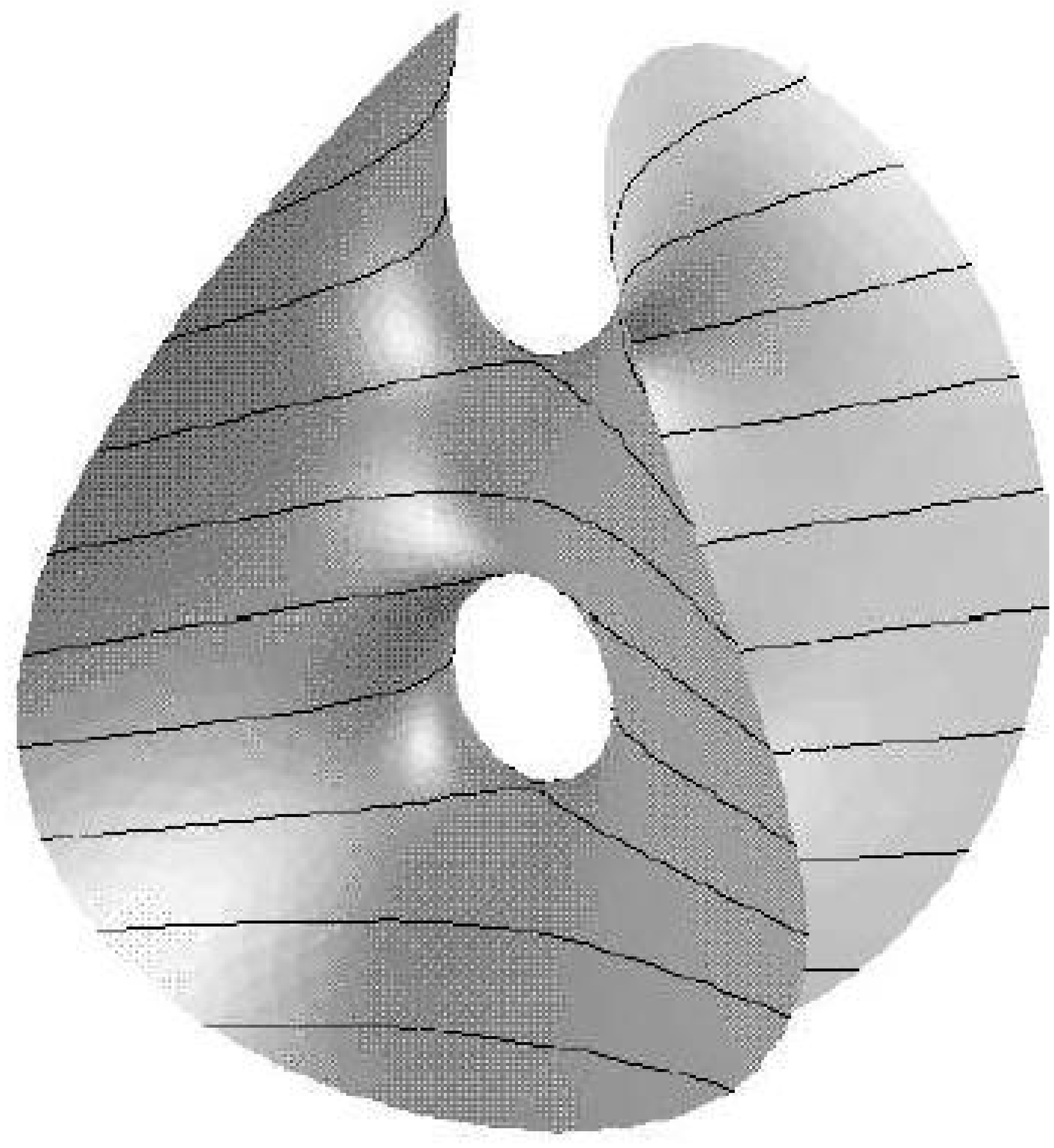,width=2.2in}
	}
}
\caption{Minimal $5$ ended surface in $\bfR^3$ with same symmetry group 
as the trinoid (proven to exist in \cite{Ka}, \cite{R}), and a minimal genus 
$1$ Enneper surface (proven to exist in \cite{cg}) in $\bfR^3$.  The 
right hand side figure was made by Edward C. Thayer of ZymoGenetics Inc. 
(Seattle, U.S.A.), and both figures were made using \cite{MESH} software.}
\label{figure7}
\end{figure}

\begin{theorem} (Weierstrass representation) 
Let $\Sigma$ be a compact Riemann surface, and $\{p_j\}_{j=1}^k \subset 
\Sigma$ be a finite number of points (representing the ends 
of the minimal surface defined in this theorem).  Fix a point 
$z_0 \in \Sigma \setminus \{p_j\}$.  
Let $g$ be a meromorphic function on $\Sigma$, and $f$ a holomorphic 
function on $\Sigma \setminus \{p_j\}$.  Assume 
that, for any point in $\Sigma \setminus \{p_j\}$, $f$ has a zero of 
order $2k$ if and only if $g$ has 
a pole of order $k$, and that $f$ has no other 
zeroes on $\Sigma \setminus \{p_j\}$.  Then, in terms of local 
holomorphic coordinates $z$ on $\Sigma \setminus \{p_j\}$, 
\begin{equation}\label{WrepR3}
 \Phi(z) = \mbox{Re} \int_{z_0}^{z}
\; \left( \begin{array}{c}
        (1-g^2)f d\zeta \\
        i(1+g^2)f d\zeta \\
        2gf d\zeta
        \end{array}
\right) \end{equation} is a 
conformal minimal immersion of the universal cover 
$\widetilde{\Sigma \setminus \{p_j\}}$ of 
$\Sigma \setminus \{p_j\}$ into $\bfR^3$. If $\Phi$ is well-defined on 
$\Sigma \setminus \{p_j\}$, then it has finite total curvature.  
Furthermore, any complete minimal surface with finite total 
curvature in $\bfR^3$ can be represented this way.
\end{theorem}

The function $g$ has a geometric interpretation: it is the composition 
of the Gauss map with stereographic projection to the complex plane $\bfC$.  
And the fundamental forms and Gaussian curvature can be described using 
$g$ and $f$: 
    \begin{equation}
     I_1 = \left( \begin{array}{cc} (1+g\bar g)^2\,f \overline{f} & 0 \\
     0 & (1+g\bar g)^2\,f \overline{f} \end{array} \right) \; , \; \; 
        I_2 = -2 \left( \begin{array}{cc} 
\mbox{Re}(\frac{dg}{du} f) & \mbox{Re}(\frac{dg}{dv} f) \\ 
\mbox{Re}(\frac{dg}{dv} f) & -\mbox{Re}(\frac{dg}{du} f) 
\end{array} \right) \; , \end{equation}\[
   K = -4 \left( \frac{|\frac{dg}{dz}|}{|f| (1+|g|^2)^2} \right)^2  \; . \]

Applying integral (1) to closed loops in $\Sigma \setminus \{p_j\}$, one 
has a period vector (a vector in $\bfR^3$) associated to each loop, 
which for a nontrivial loop might be nonzero.  This is why 
$\Phi$ might only be 
defined on a covering space of $\Sigma \setminus \{p_j\}$.  
To make a surface of finite total curvature, we need $\Phi$ well defined 
on $\Sigma \setminus \{p_j\}$ itself, which usually 
involves adjusting some real parameters in the descriptions of $f$, 
$g$ and $\Sigma \setminus \{p_j\}$ so that the real part of the integral 
in (1) about any loop in $\Sigma \setminus \{p_j\}$ is 
zero.  (If one wants periodic minimal surfaces, they can be 
constructed by having nonvanishing periods.)  

We now list six examples of Weierstrass data, and 
for the first four one can solve the integral (1) explicitly, getting 
explicit parametrizations:  
\begin{enumerate}
\item plane: $\Sigma \setminus \{p_j\}=\bfC,g=0,f=1$.  
\item Enneper surface: $\Sigma \setminus \{p_j\}=\bfC,g=z,f=1$.  
\item catenoid: $\Sigma \setminus \{p_j\}= \bfC \setminus \{0\}, 
    g=\frac{1}{z}, f = 1$.  
\item helicoid: $\Sigma \setminus \{p_j\}= \bfC \setminus \{0\}, 
    g=\frac{1}{z}, f = i$.  
\item trinoid: $\Sigma \setminus \{p_j\}=(\bfC \cup \{\infty\}) \setminus 
\{1,e^{\frac{2\pi i}{3}},e^{\frac{4\pi i}{3}}\},g=z^2,f=\frac{1}{(z^3-1)^2}$. 
\item Costa surface: $\Sigma \setminus \{p_j\}= \{(z,w) \in 
(\bfC \cup \{\infty\})^2 \; | \; w^2 = \frac{z}
{(z^2-1)} \} \setminus \{(z,w) \; | \; z = \pm 1, \infty \},g=\frac{b}{w}, 
f = \frac{z}{w},b \in \bfR$ chosen so that 
$\Phi$ is well defined on $\Sigma \setminus \{p_j\}$.  
\end{enumerate}

\section{Bryant representation}  

\begin{theorem} (Bryant representation) 
Let $\Sigma$, $\Sigma \setminus \{p_j\}$, $z$, 
$z_0$, $f$, and $g$ be as in the Weierstrass representation.  
Choose the holomorphic immersion $F:
\widetilde{\Sigma \setminus \{p_j\}} \to SL(2,\bfC)$ so 
that $F(z_0)$ is the identity matrix and $F$ satisfies 
\[ dF = F \left( \begin{array}{cc} 
g & -g^2 \\ 1 & -g \end{array} \right) f dz \; , \] then 
$\Phi:\widetilde{\Sigma \setminus \{p_j\}} \to H^3$ defined by 
\begin{equation}
        \Phi = F \cdot F^*
    \label{eq:imm}
\end{equation}
is a conformal CMC $1$ immersion in the Hermitian model of $\bfH^3$.  
Furthermore, any CMC $1$ surface with finite 
total curvature in $\bfH^3$ can be represented this way.
\end{theorem}

Following the terminology of Umehara and Yamada, $g$ is the 
{\it secondary} Gauss map of 
$\Phi$.  The fundamental forms and Gaussian curvature are again described 
using $g$ and $f$, and $I_1$ and $K$ are the same as in equation (2).  
The second fundamental form is $I_2+I_1$ (where $I_1$ and $I_2$ are 
as in equation (2)), as we saw from the Lawson correspondence.  

A significant difference from 
the Weierstrass representation is that, even when the $\Phi$ in (3) 
is well defined on $\Sigma \setminus \{p_j\}$, $g$ and $F$ might not be.  
We will soon see that the catenoid cousin is an example of this.  

Suppose $\Phi$ is indeed well defined on $\Sigma \setminus 
\{p_j\}$.  The {\em hyperbolic} Gauss map $G$ is the meromorphic function 
\[
G=\frac {dF_{11}}{dF_{21}}=\frac {dF_{12}}{dF_{22}} \; , 
\] where $F=(F_{ij})_{i,j=1,2}$.  
$G$ is single-valued by its geometric interpretation \cite{B}: 
$G$ is the image of the composition of a first map, from $z 
\in \Sigma \setminus \{p_j\}$ to the point in the 
sphere at infinity (in the Poincare model) at the opposite 
end of the oriented perpendicular geodesic ray starting from 
$\Phi(z)$ on the surface, and a second map, stereographic projection of the
sphere at infinity to the complex plane $\bfC$.  This 
geometrical property is strikingly similar to that of the $g$ in the 
Weierstrass representation.  

So the Bryant representation has two "Gauss" maps $g$ and $G$, whereas the 
Weierstrass representation had only one $g$.  One way to think about this 
is that in the Weierstrass representation, two roles are served 
by the {\em same} function $g$: the first 
is to describe the fundamental forms, and the second is to describe 
stereographic projection of the unit normal vector field.  In the Bryant 
representation, we need two {\em different} functions to serve these two 
roles: $g$ for the first, and $G$ for the second.  

Here are three simple examples of data for the Bryant 
representation:  
\begin{enumerate}
\item horosphere: $\Sigma \setminus \{p_j\}=\bfC,g=0,f=1$, like the data 
for a plane in $\bfR^3$.  
\item Enneper cousins: $\Sigma \setminus \{p_j\}=\bfC,g=z,f=k, 
k \in \bfR$, like the data for Enneper's surface in $\bfR^3$.  
$k$ is a nontrivial parameter, and the Enneper cousins for different $k$ 
do not differ by only a dilation, as would happen in the Weierstrass 
representation.  
\item catenoid cousins: $\Sigma \setminus \{p_j\} = \bfC \setminus \{0\}, 
g=z^\mu, f = 
\frac{1-\mu^2}{4 \mu} z^{-\mu-1}, \mu \in \bfR^+ \setminus \{1\}$.  $\Phi$ is 
well defined on $\bfC \setminus \{0\}$, even though $g$ (and $F$) is not.  
\end{enumerate}

The term "cousin" is often understood to mean that the CMC $1$ surface 
is locally isometric (via Lawson correspondence) to the corresponding 
minimal surface in $\bfR^3$.  The horosphere (resp. Enneper's cousin with 
$k=1$) is isometric to the plane (resp. Enneper's surface) in $\bfR^3$.  
So we could call the horosphere a "plane cousin", but it 
already has a name.  The catenoid cousin is 
locally isometric to a minimal catenoid, if one makes the coordinate 
transformation $w=z^\mu$ and an appropriate homothety of the minimal 
catenoid.  

However, this interpretation of the word "cousin" is not universal.  
For example, in the paper \cite{ruy1}, the term is used to mean what 
we will call "dual cousin" in the next section, and sometimes the term 
is also used to describe the Lawson correspondence combined with a 
conjugation by angle $\pi/2$.  

\section{Dual CMC $1$ surfaces in $\bfH^3$}

An interesting property of CMC $1$ surfaces in $\bfH^3$ is that if we change 
$F$ to $F^{-1}$ in equation (3) we get the CMC $1$ "dual" surface 
$\Phi^\# = F^{-1} (F^{-1})^*$ in $\bfH^3$.  Since 
$g=-\frac{dF_{12}}{dF_{11}}$ and $G=\frac{dF_{12}}{dF_{22}}$, changing 
$\Phi$ to $\Phi^\#$, that is, changing 
\[ F = \left( \begin{array}{cc} F_{11} & F_{12} \\ 
       F_{21} & F_{22} \end{array} \right) \longrightarrow 
 F^{-1} = \left( \begin{array}{cc} F_{22} & -F_{12} \\ 
       -F_{21} & F_{11} \end{array} \right) \] switches 
the secondary and hyperbolic Gauss maps.  (Note that one of $\Phi$ and 
$\Phi^\#$ being well defined on $\Sigma \setminus \{p_j\}$ does not imply 
the other is also.)  

Dual surfaces have useful properties, initially exploited 
by Umehara and Yamada \cite{uy5}.  One useful property is that 
period problems are easier to study for cousin dual surfaces than 
for the cousin surfaces, as we now explain:  

To produce CMC $1$ surfaces in $\bfH^3$, we could start with 
a well defined 
minimal surface in $\bfR^3$ with data $\Sigma \setminus \{p_j\}$ and 
$g$ and $f$, and insert this data into the Bryant 
representation to produce the cousin CMC $1$ surface in $\bfH^3$, and try 
to solve the period problems.  About a nontrivial 
loop $\alpha$ in $\Sigma \setminus \{p_j\}$, $F$ changes to $B_\alpha F$ for 
some constant matrix 
$B_\alpha \in SL(2,\bfC)$, so for the cousin surface to well 
defined on $\Sigma \setminus \{p_j\}$, we need 
\[ F F^* = B_\alpha F (B_\alpha F)^* = 
B_\alpha F F^* B_\alpha^* \] for all loops $\alpha$.  But this can be 
difficult to establish, given that $F$ depends on $z \in \Sigma \setminus 
\{p_j\}$, and given that $B_\alpha$ and $B_\alpha^*$ are on the "outer" part 
of the product on the right hand side.  And 
even when we can establish this, we 
produce only those special CMC $1$ surfaces where the secondary Gauss map 
is well defined on $\Sigma \setminus \{p_j\}$ -- Umehara and Yamada 
have demonstrated that such surfaces are very much a special case \cite{uy1}.  

On the other hand, let us consider the period problems of the dual cousin 
surface.  Solving the period problems requires that 
\[ F^{-1} (F^{-1})^* = (B_\alpha F)^{-1} ((B_\alpha F)^{-1})^* = 
F^{-1} B_\alpha^{-1} (B_\alpha^{-1})^* (F^{-1})^* \] for all loops 
$\alpha$.  In this case $B_\alpha$ and $B_\alpha^*$ are conveniently 
located in the "inner" part of the product on the right hand side, so we 
have a simple condition for solvability:  

\begin{lemma} ({\em $SU(2)$ condition})  
The dual surface $\Phi^\#$ is well-defined on $\Sigma 
\setminus \{p_j\}$ if 
$B_\alpha \in SU(2)$ for every closed loop $\alpha \in \Sigma 
\setminus \{p_j\}$.  
\end{lemma}

Dual cousins have another advantage.  As $g$ has 
been chosen from the Weierstrass data of a minimal surface 
well defined on $\Sigma \setminus \{p_j\}$, $g$ is also well defined on 
$\Sigma \setminus \{p_j\}$.  So once the $SU(2)$ condition is solved, 
the dual cousin's hyperbolic Gauss map, which is $g$, is well 
defined on $\Sigma \setminus \{p_j\}$, and it should be.  
But solving the $SU(2)$ 
condition does not force the secondary Gauss map, which is now $G$, to be 
well defined on $\Sigma \setminus \{p_j\}$, so the dual cousin's 
secondary Gauss map might not be well 
defined on $\Sigma \setminus \{p_j\}$, allowing for many more 
possibilities.  

We now give three examples of data for dual cousins, which are the same 
as the data for the corresponding minimal surfaces (unlike the data for the 
original cousins in non-simply-connected cases): 
\begin{enumerate}
\item Enneper cousin duals: $\Sigma \setminus \{p_j\} = 
\bfC,g=z,f=k, k \in \bfR$, like the data for Enneper's surface.  
The Enneper cousin duals, having infinite total curvature, are truly different 
from the Enneper cousins, having finite total curvature.  
\item catenoid cousin "duals": $\Sigma 
\setminus \{p_j\} = \bfC \setminus \{0\}, g=
\frac{1}{z}, f = k, k \in \bfR$, like the data for a minimal catenoid.  
The catenoid cousin "duals" are the same as the catenoid cousins 
of the previous section (hence "duals" is in quotes).  
Here we are describing the data in terms of the hyperbolic 
Gauss map, rather than the secondary Gauss map, as we did in the previous 
section.  So we see that the secondary Gauss map is not well defined.  
\item trinoid cousin duals: $\Sigma \setminus \{p_j\} 
=(\bfC \cup \{\infty\}) \setminus 
\{1,e^{\frac{2\pi i}{3}},e^{\frac{4\pi i}{3}}\}, g=z^2, f = 
\frac{k}{(z^3-1)^2}, k \in \bfR$, like the data for a trinoid.  
The secondary Gauss map of the trinoid cousin dual is also not well 
defined, and is in fact unknown.  
\end{enumerate}

\section{making dual cousins}

As many of the examples shown in the figures here have been proven 
to exist by a method in \cite{ruy1}, we briefly outline the 
method here.  For each minimal surface with finite total curvature 
satisfying certain conditions, the method implies existence of a 
one parameter family of corresponding CMC $1$ dual cousins.  

We start with a minimal surface $\Phi: \Sigma \setminus 
\{p_j\} \to \bfR^3$ of finite total curvature with Weierstrass 
data $f$ and $g$.  
We require the immersion to be symmetric in the following 
sense, a condition that generically eliminates virtually all minimal 
surfaces, but eliminates none of the better known surfaces, which 
all have symmetries:  

\begin{quote} {\bf Symmetry condition:} 
  There is a disk $D \subset \Sigma \setminus \{p_j\}$ so that 
  $\Phi(D)$ is bounded by non-straight planar geodesics.  
\end{quote}

If $\Phi$ is symmetric with respect to a disk $D$, then $\Phi(D)$ 
generates the full surface by reflections across planes containing 
the boundary planar geodesics of $\partial \Phi(D)$ 
(by the Schwartz reflection principle \cite{O}).  
Since the surface has finite total curvature, it is not
periodic, so if any two of these boundary planar geodesics lie in
parallel planes, they must lie in the same plane.  And in fact, 
it is shown in \cite{ruy1} that the boundary $\partial \Phi(D)$ is 
contained entirely in only either 
two intersecting planes $P_1$, $P_2$, or in three planes 
$P_1$, $P_2$, and $P_3$ in general position.  
Let the boundary planar geodesics of $\Phi(D)$ contained in $P_j$ 
be called $S_{j,1}$, $S_{j,2}$, \dots, $S_{j,d_j}$  
($j=1,\dots,s$, for $s=2$ or $3$).  Two examples of symmetric surfaces are: 
\begin{enumerate}
\item the genus $1$ Costa surface.  This surface can be 
placed in $\bfR^3$ so that the 
central planar end is asymptotic to the plane $\{(x_1,x_2,x_3) \in \bfR^3 
\; | \; x_3 = 0\}$ and 
so that reflections through the planes $P_1 = \{(x_1,x_2,x_3) \in \bfR^3 
\; | \; x_1=0\}$ and $P_2 = \{(x_1,x_2,x_3) \in \bfR^3 
\; | \; x_2=0\}$ are 
symmetries of the surface.  Then the piece of the surface lieing in the 
region $\{ (x_1,x_2,x_3) \in \bfR^3 \; | \; x_1,x_2 \geq 0 \}$ is 
one of four congruent disks comprising the surface.  The Costa surface 
is symmetric with respect to this disk.  The boundary of this disk is 
contained in $P_1 \cup P_2$, with 
four planar geodesics $S_{1,1}$, $S_{1,2}$, $S_{2,1}$, and $S_{2,2}$.  
\item the genus $1$ trinoid.  This surface 
can be separated into twelve congruent 
disks, any of which is bounded by four planar geodesics lying in three planes 
$P_1$, $P_2$, and $P_3$.  So the list of boundary planar geodesics of this 
disk could be written $S_{1,1}$, $S_{2,1}$, $S_{3,1}$, and $S_{3,2}$.  
\end{enumerate}

\begin{figure}
\begin{center}
\unitlength=0.8pt
\begin{picture}(300.00,220.00)(90.00,560.00)
\put(125.00,669.00){\makebox(0,0)[cc]{$\Phi(D)$}}
\put(165.00,668.00){\makebox(0,0)[cc]{\footnotesize $S_{1,1}$}}
\put(149.00,602.00){\makebox(0,0)[cc]{\footnotesize $S_{3,2}$}}
\put(90.00,592.00){\makebox(0,0)[cc]{\footnotesize $S_{2,1}$}}
\put(81.00,730.00){\makebox(0,0)[cc]{\footnotesize $S_{3,1}$}}
\bezier280(136.00,702.00)(149.00,743.00)(150.00,650.00)
\bezier160(150.00,650.00)(154.00,595.00)(160.00,620.00)
\bezier100(160.00,620.00)(130.00,621.00)(130.00,600.00)
\bezier170(130.00,600.00)(83.00,593.00)(90.00,630.00)
\bezier380(90.00,630.00)(135.00,716.00)(40.00,730.00)
\bezier20(136.00,702.00)(37.00,699.00)(40.00,730.00)
\put(360.00,780.00){\line(0,-1){160.00}}
\put(360.00,620.00){\line(-1,0){0.00}}
\put(290.00,740.00){\line(0,-1){160.00}}
\put(290.00,580.00){\line(5,3){70.00}}
\put(360.00,622.00){\line(1,-2){31.00}}
\put(391.00,560.00){\line(0,1){160.00}}
\put(391.00,720.00){\line(-1,2){30.00}}
\put(361.00,780.00){\line(-5,-3){71.00}}
\put(341.00,669.00){\makebox(0,0)[cc]{$\Phi(D)$}}
\put(377.00,712.00){\makebox(0,0)[cc]{\footnotesize $S_{1,1}$}}
\put(388.00,615.00){\makebox(0,0)[cc]{\footnotesize $S_{1,2}$}}
\put(318.00,723.00){\makebox(0,0)[cc]{\footnotesize $S_{2,1}$}}
\put(327.00,620.00){\makebox(0,0)[cc]{\footnotesize $S_{2,2}$}}
\bezier300(390.00,694.00)(356.00,732.00)(390.00,635.00)
\bezier380(290.00,661.00)(375.00,666.00)(290.00,605.00)
\bezier110(360.00,685.00)(373.00,662.00)(374.00,635.00)
\bezier120(374.00,635.00)(372.00,616.00)(391.00,582.00)
\bezier120(361.00,685.00)(334.00,675.00)(322.00,705.00)
\bezier90(322.00,705.00)(321.00,719.00)(290.00,719.00)
\bezier20(390.00,694.00)(300.00,686.00)(290.00,719.00)
\bezier20(390.00,635.00)(300.00,626.00)(290.00,661.00)
\bezier20(391.00,582.00)(300.00,570.00)(290.00,605.00)
\end{picture}
\end{center}
    \caption{The disks of symmetry for the trinoid and Costa surface}
    \label{fig:costafund}
\end{figure}
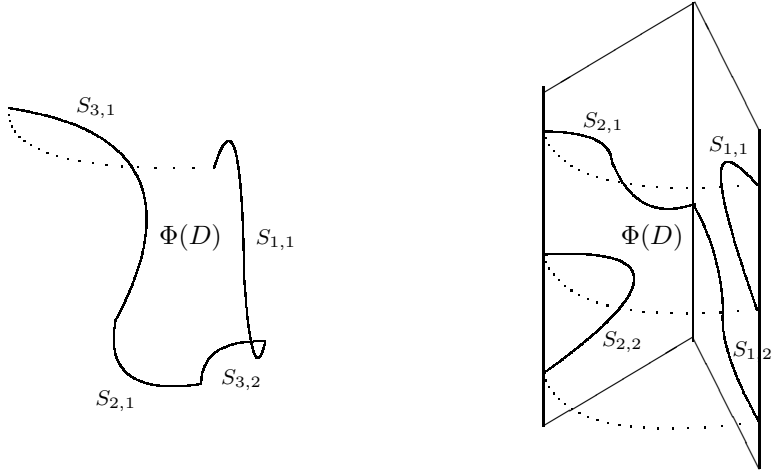

We now define non-degeneracy of the period problems.  
Let $d$ be the number of $S_{j,\ell}$ minus the number of planes 
($d= d_1+d_2+d_3-3$ if $s=3$, and $d = d_1+d_2-2$ if $s=2$).  

\begin{quote} {\bf Nondegeneracy condition:} 
There exists a continuous $d$-parameter family of 
minimal disks $\Phi(D)_\lambda$, 
  ($\lambda=(\lambda_1,...,\lambda_{d}), \, \lambda \approx \vec{0}$) such 
  that
  \begin{enumerate}
    \item   $\Phi(D)_{(0,0,...,0)} = \Phi(D)$.
    \item   $\partial \Phi(D)_\lambda = \cup_{j=1}^{s}(\cup_{\ell=1}^{d_j}
              S_{j,\ell}(\lambda))$  
            and each $S_{j,\ell}(\lambda)$ is 
            a planar geodesic lying in a plane
            $P_{j,\ell}(\lambda)$ parallel to $P_j$.  
    \item   letting $Per_{j,\ell}(\lambda)$ 
            ($j=1,\dots,s$, $\ell=2,\dots,d_j$) be the oriented distance 
            between the plane $P_{j,\ell}(\lambda)$ and $P_{j,1}(\lambda)$, 
            the map from $\lambda$ in $\bfR^d$ to 
            $(Per_{j,\ell}(\lambda))$ in $\bfR^d$ is 
            an open map onto a neighborhood of $\vec{0}$.  
  \end{enumerate}
\end{quote}

\begin{theorem}\label{mainthm} \cite{ruy1} 
If the minimal immersion $\Phi$ is symmetric and 
nondegenerate, then there exists an associated 
one-parameter family of CMC $1$ cousin duals in $\bfH^3$.  
\end{theorem}

\begin{quote}
{\small 

We outline the proof:  We will not show that the period problems 
are solvable, but we will at least 
explain why we have enough free parameters available so that the 
period problems are not overdetermined.  

\begin{figure}
\centerline{
        \hbox{
		\psfig{figure=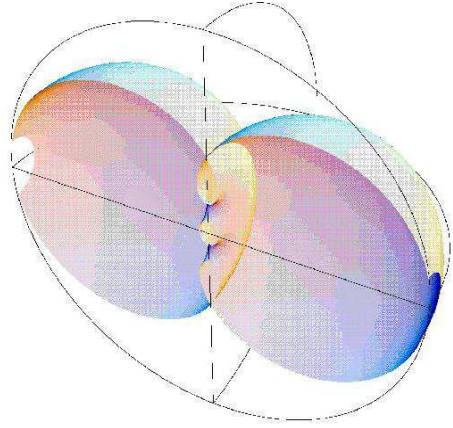,width=2.5in}
	}
}
\caption{
CMC $1$ genus $1$ catenoid cousin dual in $\bfH^3$, proven to exist in 
\cite{rs}.  The graphics were made 
by Katsunori Sato of Tokyo Institute of Technology.  
No corresponding minimal surface can exist, by Schoen's result (which is 
mentioned in the caption of Figure 8).  
Levitt and Rosenberg \cite{LR} have proved that any complete properly 
embedded CMC $1$ surface in $\bfH^3$ with asymptotic boundary 
consisting of at most two points is a surface of revolution, which implies 
that this example, and the examples in Figures 6 and 7, cannot 
be embedded, and we see that they are not.}
\label{figure18}
\end{figure}

\begin{figure}
\vspace{0.5in}
\centerline{
        \hbox{
		\psfig{figure=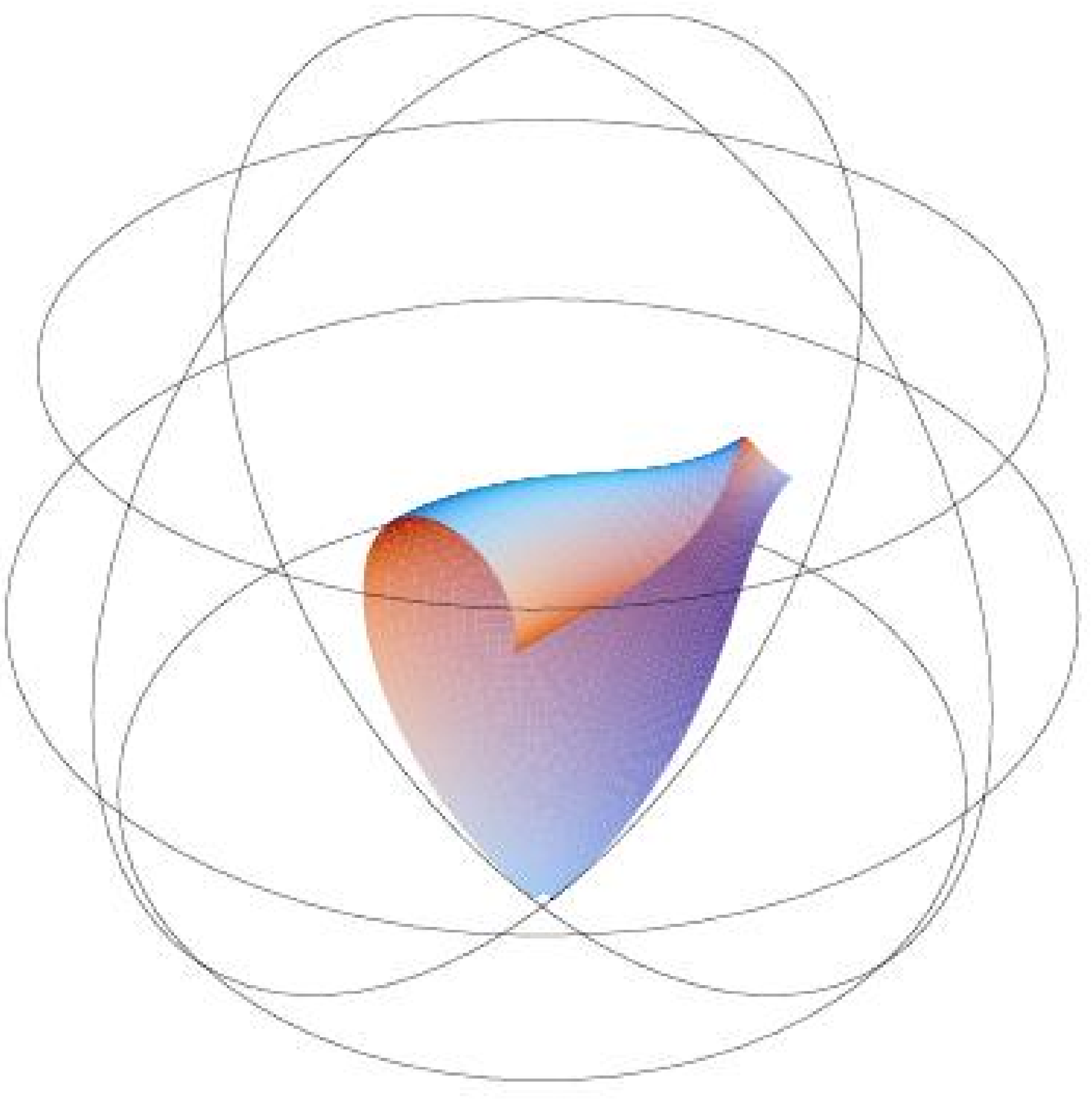,width=2in}
		\hspace{0.6in}
		\psfig{figure=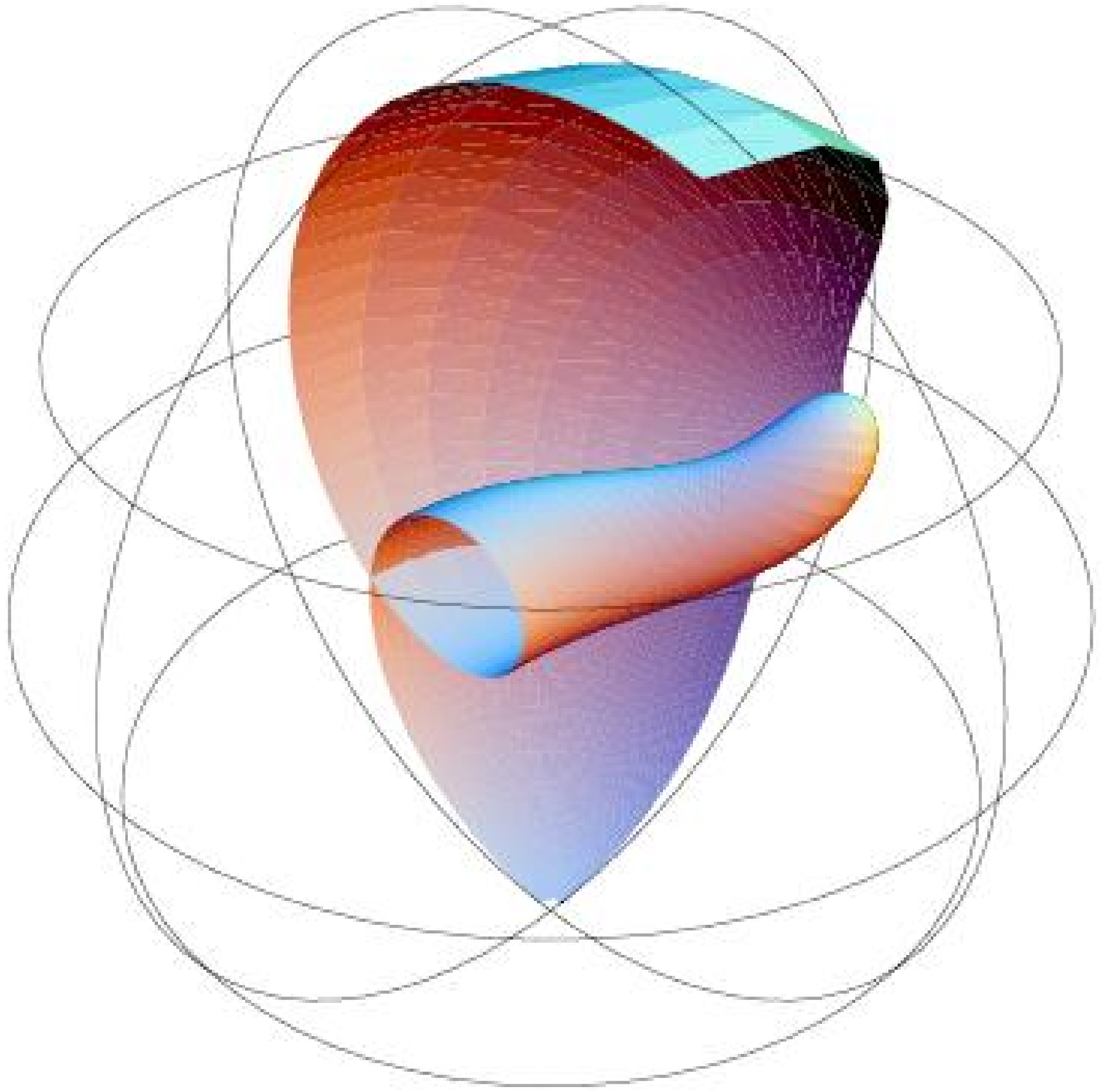,width=2in}}}
		\vspace{0.5in}
		\centerline{\hbox{
		\psfig{figure=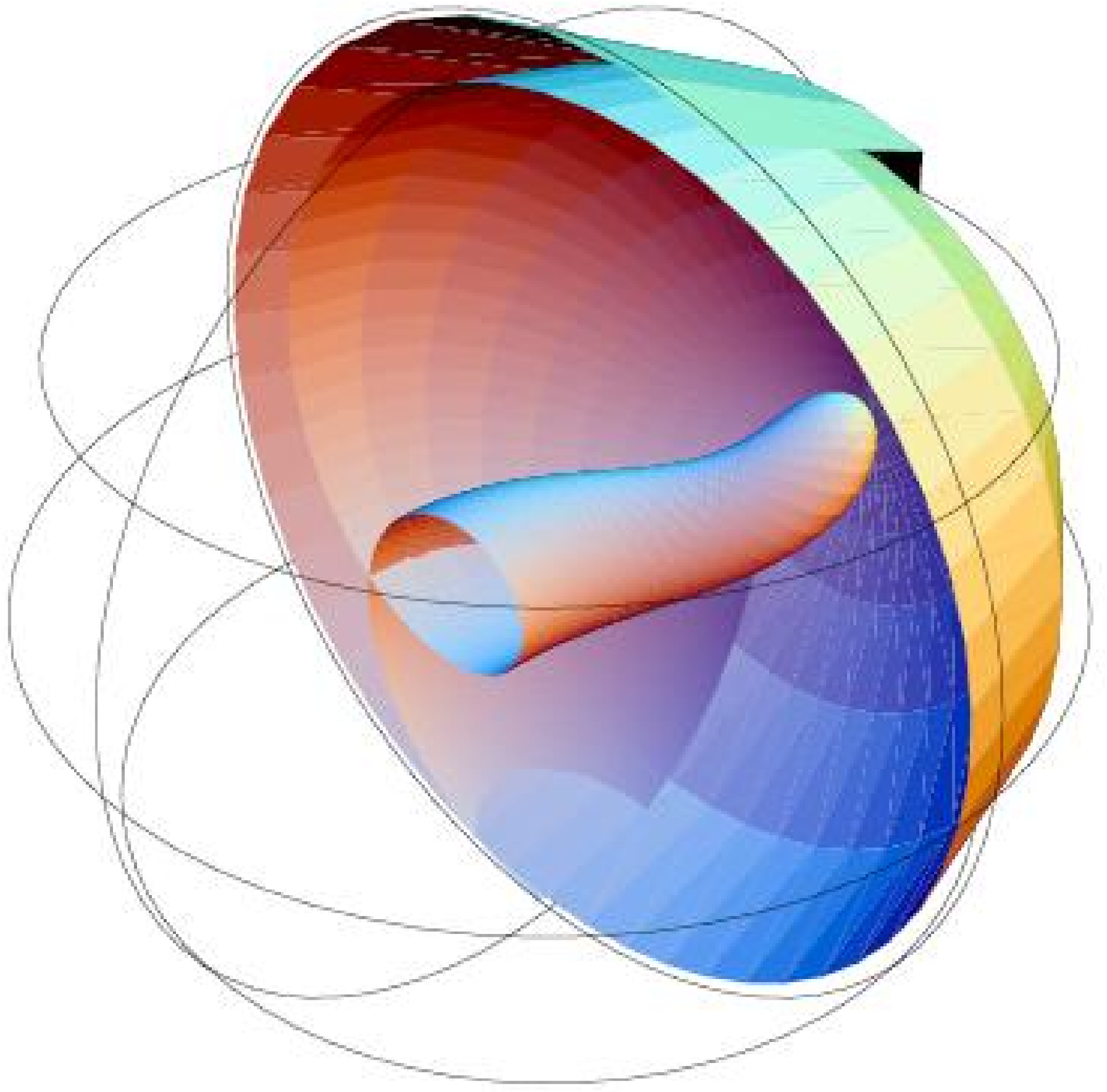,width=2in}
	}
}
\caption{
CMC $1$ surface in $\bfH^3$, proven to exist in \cite{uy1}.  This example 
is interesting because the hyperbolic Gauss map has an essential singularity 
at one of its two ends, like the end of the Enneper cousin.  And the 
geometric behavior of the end here is strikingly similar to that 
of the Enneper cousin's end.  Here we show three pictures consecutively 
including more of this end.}
\label{figure19}
\end{figure}

First note that studying 
CMC $1$ surfaces in the standard hyperbolic space 
is equivalent to studying CMC $c$ surfaces in the hyperbolic 
space of constant sectional curvature $-c^2$.  (One can see this by 
dilating $\bfR^3$ from the origin by a factor 
$\frac{1}{c}$, taking the standard Poincare model to the model 
$\{(x_1,x_2,x_3) \in \bfR^3 \; | \; x_1^2+x_2^2+x_3^2 < \frac{1}{c} \}$ with 
metric $ds_c^2 = \frac{4 \sum dx_i^2}{(1-c^2 \sum x_i^2)^2}$, which is 
hyperbolic space of constant sectional curvature $-c^2$, and 
taking CMC $1$ surfaces to CMC $c$ surfaces.)  So we can instead study 
the problem with $c \approx 0$, where we are looking 
at only slight deformations of the original minimal surface, so that 
the nondegeneracy of the minimal surface's periods can be used.  

Now suppose $D \subset \Sigma \setminus \{p_j\}$ is a disk by which the 
minimal surface is symmetric.  
Consider a local coordinate $\{z \in \bfC \; | \; |z|<1 \}$ of 
$\Sigma$ in a neighborhood of 
a point in $\partial D$, which can be chosen so that $\{z \in \bfR \; | \; 
|z|<1 \}$ lies on $\partial D$ (hence $\Phi$ maps it to a portion of 
some $S_{j,\ell} \subset P_j \cap \partial \Phi(D)$), and so that 
$z \to \bar{z}$ represents reflection of the surface across the plane 
$P_j$.  It turns out that the $F$ from the Bryant representation satisfies 
\[ \overline{F(\bar{z})} = \rho_{j,\ell} F(z) \sigma_j \; , \]\[ 
               \sigma_j =  \begin{pmatrix}  
            n_{2,j} - i n_{1,j} & i n_{3,j} \\
            i n_{3,j}  & n_{2,j}+i n_{1,j}
          \end{pmatrix}
          \; , \; \; 
          \rho_{j,\ell} = \begin{pmatrix}
                        p_{j,\ell} &   i\gamma_{1,j,\ell}\\
                        i\gamma_{2,j,\ell}  &   \overline{p_{j,\ell}}
                     \end{pmatrix}
          \; , 
           \]                      
where $(n_{1,j},n_{2,j},n_{3,j})$ is a unit normal vector 
to the plane $P_j$, and $\sigma_j, \rho_{j,\ell} \in 
SL(2,\bfC)$, and where $p_{j,\ell} \in \bfC$ and $\gamma_{1,j,\ell}, 
\gamma_{2,j,\ell} \in \bfR$.  One can now think of loops in 
$\Sigma \setminus \{p_j\}$ as 
compositions of reflections (a reflection for each $S_{j,\ell}$), and, 
with a bit of consideration regarding how $F$ would change under the 
compositions of such reflections that represent closed loops, we can 
see that the matrices $\sigma_j$ never create any 
period problems, and that the $SU(2)$ condition is satisfied if each 
$\rho_{j,\ell} \in SU(2)$, i.e. 
$\gamma_{1,j,\ell} = \gamma_{2,j,\ell}$ for all $S_{j,\ell}$.  Hence the 
period problem 
is $1$ dimensional for each $S_{j,\ell}$, and so in full it 
is a $d+s$ dimensional problem, as there are $d+s$ curves $S_{j,\ell}$.  

But, as we saw in the definition of nondegeneracy of the original minimal 
surface, there we only had a 
$d$ real dimensional period problem, so the Weierstrass 
data of the original minimal surface only guarantees us $d$ free 
parameters $\lambda_1,...,\lambda_d$.  We need an extra 
$s$ free parameters to prevent 
the period problems from being overdetermined, and we get them by changing 
the initial condition $F(z_0)$ in the Bryant representation from the 
identity to other matrices in $SL(2,\bfC)$.  (Changing the initial condition 
nontrivially changes the resulting dual cousin in $\bfH^3$, even though 
a change in the initial condition of the Weierstrass 
representation amounts to only a trivial translation in $\bfR^3$, and 
a change in the initial condition of the Bryant representation 
for the original cousins amounts to only a trivial isometric motion of 
$\bfH^3$.)  The upshot is that we have $d+s$ free 
parameters available for each value of $c$, and, as shown in \cite{ruy1}, 
the $SU(2)$ condition is solvable for each $c \approx 0$, using the 
nondegeneracy condition of the minimal surface.  So we have 
a $1$ parameter family of dual cousins in $\bfH^3$ with parameter 
$c \approx 0$.  

}
\end{quote}

Since many types of minimal surfaces in $\bfR^3$ are symmetric and 
nondegenerate, this result implies the existence of many 
types of CMC $1$ dual cousins, some of which are shown in the figures.  

\section{Other aspects of CMC $1$ surfaces in $\bfH^3$}

Up to this point, we have been considering how to make 
examples of CMC $1$ immersions in $\bfH^3$, a primary 
interest of the author.  But of course there are other 
approaches to the study of these surfaces, so we close with a brief 
list of other research and viewpoints that 
the author is aware of.  Some terms in this section we leave undefined, 
but the reader can refer to the original sources.  

First we note that nice general introductions to the subject can be found in 
\cite{uy6} and \cite{Ro}.  In \cite{uy6}, like here, comparison with minimal 
surfaces in $\bfR^3$ is used, but the presentation and emphasized points 
are quite different than in this paper.  In \cite{Ro}, there is a detailed 
explanation of the Bryant representation and of some fundamental examples.  

We begin with some of the results by Umehara and Yamada -- who have been 
pursuing a fruitful long-term study of CMC $1$ surfaces in $\bfH^3$ -- as 
their results bear the closest relation with the previous sections.  

\begin{itemize}

\item In \cite{uy1}, Umehara and Yamada use the Frobenius method to prove a 
number of results.  Among the results in \cite{uy1} are the following: 
\begin{enumerate}
\item reducible CMC $1$ surfaces in $\bfH^3$ have a Perez-Ros type deformation 
\cite{PR}, i.e. changing $g,f$ to $\lambda g,f/\lambda$ for $\lambda \in 
\bfR \setminus \{ 0 \}$ leaves $\Phi$ still well-defined on 
$\Sigma \setminus \{ p_j \}$; 
\item for CMC $1$ surfaces in $\bfH^3$, equality 
never holds in the Cohn-Vossen inequality 
\[ \frac{1}{2\pi} \int_\Sigma |K| \ dA > -\chi (\Sigma \setminus \{ p_j \} ) ; 
\]
\item there is a simple algebraic way to check if a complete regular end is 
embedded, relying only on the orders and the leading coefficients of 
$f$ and $g$ at the end; 
\item all genus 0 immersions with 2 regular ends are 
classified; 
\item a general procedure for solving period problems at regular ends is 
given, which is applied to make CMC $1$ genus 0 surfaces; 
\item any regular end is shown to be tangent to the sphere at infinity at a 
single limiting point.  
\end{enumerate}

\item In \cite{uy3}, a correspondence between complete CMC $1$ surfaces in 
$\bfH^3$ with finite total curvature and abstract surfaces with constant 
Gauss curvature $1$ and isolated conical singularities (Met$_1$ surfaces) 
is established, and 
also CMC $1$ surfaces with dihedral symmetry (similar to the Jorge-Meeks 
minimal surfaces in $\bfR^3$ \cite{JM}) and 
with the same symmetry as the Platonic 
solids (similar to surfaces in \cite{Xu} and \cite{Ka}) are shown to exist.  

\item In \cite{uy5} and \cite{Yu2}, dual 
CMC $1$ surfaces are explained (these are the 
sources for Section 5 of this paper).  
Using duality, these papers establish an Osserman type 
inequality 
\[ \frac{1}{2\pi} \int_\Sigma |K^\#| \ dA^\# \geq 
-\chi (\Sigma \setminus \{ p_j \} ) + 
\mbox{number of ends of original surface} 
\]
for the total curvature of the dual surfaces, analogous to the 
Osserman inequality for minimal surfaces in $\bfR^3$.  ($K^\#$ and $dA^\#$ are 
the Gaussian curvature and area form of the dual surfaces.)  A 
point of interest here is that the Osserman inequality does not hold for the 
total curvature $\int_\Sigma K dA$ of the original CMC $1$ surfaces in 
$\bfH^3$ (the Cohn-Vossen inequality is the best possible result), 
and one really 
must go to the dual surfaces to establish the Osserman inequality. 

\item In \cite{uy7}, the correlation between CMC $1$ surfaces in 
$\bfH^3$ and abstract Met$_1$ surfaces is used to classify 
Met$_1$ surfaces with $3$ conical singularities.  Also, genus 
$0$ CMC $1$ irreducible immersions with $3$ embedded 
ends all of which are asymptotic 
to catenoid cousin ends ("trinoids") are 
classified.  We remark that Collin and Rosenberg 
have recently been creating a more geometric method for classifying trinoids. 

\item With the author, Umehara and Yamada have recently used residue 
conditions of the data in the dual surfaces' Bryant representation 
to show nonexistence 
of certain types of genus $0$ CMC $1$ surfaces in $\bfH^3$ \cite{ruy2}.  

\item Umehara and Yamada and the author have been 
classifying CMC $1$ surfaces with low total curvature, 
and surfaces whose dual surfaces have low
total curvature \cite{ruy3}, \cite{ruy4}, \cite{ruy5}, \cite{ruy6}.  
(The Osserman inequality gives one reason for the interest in surfaces whose 
dual surfaces have low total curvature.)  

\end{itemize}

The flux of CMC surfaces in $\bfH^3$ is defined and explored by Korevaar, 
Kusner, Meeks and Solomon in \cite{kkms}, and in particular this gives a 
flux when $H \equiv 1$.  
The residue conditions in \cite{ruy2} 
have similar properties to the flux in \cite{kkms} when $H \equiv 1$, 
making it believable that the residue conditions 
are equivalent to the flux, but this is still unproven.  

An interesting result by Small \cite{Sm}, using algebraic techniques, 
shows that when $g$ and $G$ are known, the surface can be given explicitly 
in terms of $g$, $G$, $G^\prime$, and $G^{\prime\prime}$.  (Small's result is 
stated in \cite{uy6} in the notation of Umehara and Yamada.)  

McCune and Umehara \cite{MU} have recently found an analogy for CMC 
$1$ surfaces in $\bfH^3$ of the UP iteration for making minimal surfaces in 
$\bfR^3$ \cite{M}.  With it, they make new examples of CMC $1$ immersions 
in $\bfH^3$.  

The results above are about immersions.  Other 
research has, however, focused on properties that CMC 
$1$ {\em embeddings} must have, and we now list some results 
about this.  We begin with an important recent work \cite{chr1}.  

\begin{itemize}
\item Recently Collin, Hauswirth, and Rosenberg have given 
an impressive proof that any properly embedded CMC $1$ surface in $\bfH^3$ of 
finite topology has finite total curvature, and all of its ends are regular 
\cite{chr1}, \cite{Ro}.  
(An end $p_j$ is regular when $G$ is meromorphic at $p_j$, or equivalently 
order$_{p_j}($Hopf differential$=g^\prime f dz^2) \geq -2$ \cite{B}.)  
Then the result of Earp and Toubiana \cite{ET1} (see below) implies that each 
end is asymptotic 
to a catenoid cousin end or a horosphere end.  (Finite total curvature and 
regularity are precisely the needed conditions in \cite{ET1}.  Also, we note 
here that this asymptotic behavior is implicit in the 
computations in Section 5 of \cite{uy1}, as argued in the appendix of 
\cite{LiRo}.) Among the corollaries of the \cite{chr1} result are:
\begin{enumerate}
\item a properly embedded CMC $1$ surface in $\bfH^3$ of 
finite topology that is simply-connected (resp. annular) must be a 
horosphere (resp. catenoid cousin); 
\item irregular CMC $1$ ends cannot be embedded, confirming a conjecture in 
\cite{uy5}.  
\end{enumerate}

Also in \cite{chr1}, it is shown that a properly embedded CMC $1$ surface in 
$\bfH^3$ of finite topology with an end asymptotic to a horosphere end must be 
an actual horosphere.  (This shows that many types of CMC 
$1$ immersions cannot be embedded, and contrasts with the fact 
that properly embedded nonplanar minimal surfaces of finite topology 
in $\bfR^3$ can have planar ends.)  

\item Given a CMC $1$ embedding $\Phi$ into $\bfH^3$ in the Poincare model 
$B^3$, let 
$\partial_\infty \Phi$ denote the limit set of $\Phi$ in the 
sphere at infinity $\partial_\infty B^3$ of 
the Poincare model.  Using the maximum principle, 
Levitt and Rosenberg \cite{LR} showed that if $\partial_\infty \Phi$ lies 
in a single great circle of $\partial_\infty B^3$, then the surface has 
a reflective symmetry with respect to the geodesic plane whose limit set 
at infinity is that great circle.  In other words, in this case the surface 
inherits the symmetry of its boundary.  (In \cite{LR}, the result is stated 
more generally for CMC $H$ hypersurfaces in $\bfH^n$.) 

\item Do Carmo and Lawson showed that a complete CMC hypersurface which is a 
proper embedding $\Phi$ into $\bfH^n$ with 
exactly one point in $\partial_\infty \Phi$ must be a horosphere.  They 
assume only that $H$ is constant, and not necessarily that it is $1$.  They 
also give other results for the $H \neq 1$ case -- for example, if the 
hypersurface $\Phi$ is 
compact, then $\Phi$ is a sphere (known by Alexandrov); 
and if $\partial_\infty \Phi$ is a hypersphere in $\partial_\infty 
\bfH^n$ and $\Phi$ separates the two components of $\partial_\infty \bfH^n 
\setminus \partial_\infty \Phi$, then $\Phi$ is a hypersphere.  

\item Rodriguez and Rosenberg \cite{RR} showed that a properly embedded CMC 
$1$ surface $\Phi$ in $\bfH^3$ lieing in the interior region bounded by a 
horosphere is itself a horosphere.  Furthermore, if $\Phi$ lies in the 
exterior region bounded by a horosphere, with mean curvature vector pointing 
toward that horosphere, then $\Phi$ again is a horosphere.  

\item Earp and Toubiana \cite{ET1} showed that an embedded CMC 
$1$ end that is regular ($G$ extends meromorphically to the end) and of 
finite total curvature must be asymptotic to either 
a catenoid cousin end or a horosphere end.  They also 
include an explicit description of CMC $1$ helicoids in $\bfH^3$.  (These 
helicoids are also described in \cite{Ro}.)  

\item Do Carmo, Gomes and Thorbergsson \cite{CGT} considered 
complete properly embedded CMC hypersurfaces $\Phi$ in 
$\bfH^{n}$.  When $H \in [0,1)$, they showed that $\partial_\infty \Phi$ has 
no isolated points.  They also showed that 
for any H, if $\partial_\infty \Phi$ is $C^2$ regular at infinity, then 
$H > 1$ if and only if $\partial_\infty \Phi$ is empty.  (Ends that are 
asymptotic to hyperbolic Delaunay surfaces, as in \cite{kkms}, are not 
$C^2$ regular at infinity.)  These results 
contrast with the case of $H \equiv 1$, where $\partial_\infty \Phi$ is 
never empty and can have isolated points.  Another result is that if 
$H > 1$, $\partial_{\infty} \Phi$ does not contain any component of 
codimension $1$ in $\partial_\infty \bfH^n$.  

\item Karcher \cite{K} has recently constructed several different types of new 
embedded CMC $1$ surfaces in $\bfH^3$ of infinite topology, using properties 
of conjugate minimal surfaces in $\bfR^3$.  These surfaces have the same 
symmetries as Platonic tessellations of $\bfH^3$.  

\end{itemize}

Another line of inquiry is the stability and Morse index of 
CMC $1$ surfaces in $\bfH^3$.  We saw in Section 1 that area of a CMC 
immersion $\Phi:\Sigma \setminus \{ p_j \} \to \bfH^3$ is 
critical for any volume-preserving compactly supported variation 
$\Phi_t$ of $\Phi$.  So to determine if $\Phi_t$ reduces area, we must 
use the second variation formula.  
By reparametrizing the surfaces of $\Phi_t$ if 
necessary, we may assume that 
\[ \left. (\Phi_t)_* \frac{\partial}{\partial 
t} \right|_{t=0} = u \vec{N} \] for some compactly supported smooth function 
$u$ on $\Sigma \setminus \{ p_j \}$ and 
that $\int_{\Sigma} u dA = 0$, and then the second 
variation formula becomes 
\[ \left. \frac{d^2}{dt^2} \mbox{Area}(\Phi_t(U)) \right|_{t=0} = 
\int_\Sigma u \cdot (\triangle + 2 K)u \ dA \; . \]  If this integral is 
nonnegative for all compactly supported $u$ with (resp. without) the condition 
$\int_\Sigma u dA = 0$, then the surface is called {\em stable} (resp. 
{\em strongly stable}).  The {\em index} of the surface is the maximum 
dimension of a vector space of compactly supported smooth functions $u$ such 
that the second variation formula 
is negative and $\int_\Sigma u dA = 0$, for all nonzero functions $u$ in the 
space.  (Hence the surface is stable if and only if the index is $0$.)  

\begin{itemize}
\item Silveira \cite{Si} showed that for a complete noncompact CMC immersion 
in $\bfH^3$ with $H \geq 1$, stability of the surface implies that $H=1$ and 
the surface is a horosphere.  He also showed that a compact CMC $1$ disk with 
boundary and total curvature less than $2\pi$ must be 
strongly stable, analogous to the result of \cite{BC} for minimal 
surfaces in $\bfR^3$.  

\item Do Carmo and Silveira \cite{CS} showed that a CMC $1$ surface in 
$\bfH^3$ has finite index if and only if it has finite total curvature, 
analogous to the same result for minimal surfaces in $\bfR^3$ by 
Fischer-Colbrie \cite{FC} and Gulliver \cite{G}.  Lima and the author 
\cite{LiRo} followed this with index estimates for 
specific examples, in particular the index of the catenoid cousins is 
computed.  

\item There is also a concept of {\em strong index} for CMC $1$ surfaces in 
$\bfH^3$, which is the number of negative eigenvalues of the operator 
$\triangle + 2 K$ in the second variation formula with 
an appropriate domain of 
functions (strong index is explained in \cite{CS}, \cite{LiRo}, and 
\cite{BB}).  Strong index is always at least as big as the index, 
but is at most one greater when the index is finite.  So the strong index 
is a good estimate for the index, and has the advantage of being easier 
to compute, with its more analytic nature.  So it is interesting to know when 
index and strong index agree, and in fact, it is shown 
by Barbosa, Berard, and Hauswirth \cite{BB}, \cite{BH} that they always 
agree for complete CMC $1$ surfaces in $\bfH^3$.  

\end{itemize}

We remark that Goes, Galvao and Nelli 
\cite{ggn}, and also Earp and Toubiana \cite{ET2}, have developed 
alternatives to the Bryant representation, used to find examples.  

Finally we note that there are results on the omitted points of the image of 
the hyperbolic Gauss map $G$.  Yu \cite{Yu1} has shown that 
for a complete CMC $1$ surface which is not a horosphere, $G$ omits at most 
$4$ points.  And recently Collin, Hauswirth and Rosenberg \cite{chr2} have 
shown that a complete CMC $1$ immersion with finite total curvature 
that is not a horosphere has $G$ missing 
at most $3$ points.  Also, in \cite{chr2} it is shown that 
a properly embedded CMC $1$ surface of finite topology that is not a 
horosphere or catenoid cousin has surjective $G$.

wayne@math.kobe-u.ac.jp, www.math.kobe-u.ac.jp/HOME/wayne/wayne.html

\end{document}